\documentstyle[amssymb,amstex]{article}
\newtheorem{theorem}{Theorem}
\newtheorem{lemma}{Lemma}
\newtheorem{corollary}{Corollary}
\newtheorem{proposition}{Proposition}
\newcommand{\lgth}{\operatorname{length}}
\newcommand{\card}{\operatorname{card}}

\newcommand{\iin}{\operatorname{in}}
\newcommand{\oout}{\operatorname{out}}
\newcommand{\innt}{\operatorname{Int}}
\newcommand{\Var}{\operatorname{Var}}
\newcommand{\Po}{\operatorname{Po}}
\newcommand{\dist}{\operatorname{dist}}
\newcommand{\conv}{\operatorname{conv}}
\newcommand{\diam}{\operatorname{diam}}
\newcommand{\Ver}{\operatorname{Vertices}}
\newcommand{\Cov}{\operatorname{Cov}}
\newcommand{\e}{\operatorname{e}}
\begin{document}
\setlength{\baselineskip}{1.2\baselineskip}

\title{Random dynamics and thermodynamic limits for polygonal Markov fields in the
 plane}
\author{Tomasz Schreiber\footnote{Research supported by the Foundation for Polish
         Science (FNP)},\\
        Faculty of Mathematics and Computer Science,\\
        Nicolaus Copernicus University,\\
        Toru\'n, Poland,\\
        {\sc e-mail}: tomeks at mat.uni.torun.pl}
\date{}
\maketitle

\paragraph{Abstract:}
 {\it We construct random dynamics on collections of non-intersecting
      planar contours, leaving invariant the distributions of length-
      and area-interacting polygonal Markov fields with V-shaped nodes.
      The first of these dynamics is based on the dynamic construction
      of consistent polygonal fields, as presented in the original
      articles by Arak (1982) and Arak \& Surgailis (1989, 1991), and
      it provides an easy-to-implement Metropolis-type simulation algorithm.
      The second dynamics leads to a graphical construction in the spirit of
      Fern\'andez, Ferrari \& Garcia (1998,2002) and it yields a perfect
      simulation scheme in a finite window from the infinite-volume limit.
      This algorithm seems difficult to implement, yet its value lies in
      that it allows for theoretical analysis of thermodynamic limit
      behaviour of length-interacting polygonal fields. The results
      thus obtained include the uniqueness and exponential
      $\alpha$-mixing of the thermodynamic limit of such fields in
      the low temperature region, in the class of infinite-volume
      Gibbs measures without infinite contours. Outside this class
      we conjecture the existence of an infinite number of extreme
      phases breaking both the translational and rotational
      symmetries.}
\paragraph{Keywords:}
 {\it Polygonal Markov fields, random dynamics, Metropolis simulation, perfect
      simulation, thermodynamic limit, phase transitions}
\paragraph{MSC: 60D05, 60K35, 82B21}

\section{Introduction}
 An example of a planar Markov field with polygonal realisations was
 first introduced in Arak (1982). The original Arak process in a bounded
 open convex set $D$ is constructed as briefly sketched below. We define
 the family $\Gamma_D$ of admissible polygonal configurations on $D$
 by taking all the finite planar graphs $\gamma$ in $D \cup \partial D,$
 with straight-line segments as edges, such that
 \begin{description}
  \item{\bf (P1)} the edges of $\gamma$ do not intersect,
  \item{\bf (P2)} all the interior vertices of $\gamma$ (lying in $D$) are of degree $2,$
  \item{\bf (P3)} all the boundary vertices of $\gamma$ (lying in $\partial D$) are of degree $1,$
  \item{\bf (P4)} no two edges of $\gamma$ are colinear.
 \end{description}
 In other words, $\gamma$ consists of a finite number of disjoint polygons,
 possibly nested and chopped off by the boundary. Further, for a finite collection
 $(l) = (l_i)_{i=1}^n$ of straight lines intersecting $D$ we write $\Gamma_D(l)$ for the
 family of admissible configurations $\gamma$ with the additional properties that
 $\gamma \subseteq \bigcup_{i=1}^n l_i$ and $\gamma \cap l_i$ is a single interval
 of a strictly positive length for each $l_i, i=1,...,n,$ possibly with some
 isolated points added. Let $\Lambda_D$ be the
 restriction to $D$ of a homogeneous Poisson line process $\Lambda$ with intensity measure
 given by the standard isometry-invariant Lebesgue measure $\mu$ on the space of
 straight lines in ${\Bbb R}^2.$ One possible construction of $\mu$ goes by
 identifying a straight line $l$ with the pair $(\phi,\rho) \in [0,\pi) \times {\Bbb R},$
 where $(\rho \sin(\phi), \rho \cos(\phi))$ is the vector orthogonal to $l,$ and joining
 it to the origin, and then by endowing the parameter space
 $[0,\pi) \times {\Bbb R}$ with the usual Lebesgue measure. With the above notation,
 the polygonal Arak process ${\cal A}_D$ on $D$ arises as the Gibbsian modification
 of the process induced on $\Gamma_D$ by $\Lambda_D,$ with the Hamiltonian given
 by the double total edge length, that is to say  
 \begin{equation}\label{GREPR1}
    {\Bbb P}({\cal A}_D \in G)  = \frac{{\Bbb E} \sum_{\gamma \in \Gamma_D(\Lambda_D) \cap G}
     \exp(-2\lgth(\gamma))}{{\Bbb E}\sum_{\gamma \in \Gamma_D(\Lambda_D)}
                                      \exp(-2 \lgth(\gamma))}
 \end{equation}
  for all $G \subseteq \Gamma_D$ Borel measurable, say with respect to the usual
  Hausdorff distance topology, see Section 4 in Arak \& Surgailis (1989).
 The Arak process has a number of remarkable properties. It is exactly solvable
 (an explicit formula for the partition function is available), consistent
 (${\cal A}_D$ coincides in distribution with the restriction of ${\cal A}_C$
 to $D$ for $C \supseteq D$) and enjoys a two-dimensional Markov property stating
 that the conditional behaviour of the process in an open bounded domain depends on
 the exterior configuration only through arbitrarily close neighbourhoods of the
 boundary, see ibidem.
 These nice features are shared by a much broader class of processes, so-called
 consistent polygonal Markov fields, introduced and investigated in detail in Arak \&
 Surgailis (1989, 1991). Arak, Clifford \& Surgailis (1993) introduce an alternative
 point- rather than line-based representation of these models. Our description below
 specialises for the standard Arak process ${\cal A}_D.$ For a given point configuration
 $\bar{x} = \{ x_1,...,x_n \} \subseteq D \cup \partial D$ denote by $\Gamma_D(\bar{x})$
 the family of admissible configurations $\gamma$ whose vertex set coincides with
 $\bar{x}.$ Write $\Pi_D$ for the Poisson point process in $D \cup \partial D$ with
 the intensity measure given by the area element on $D$ and by the length element
 on $\partial D.$ By Theorem 1 ibidem (see also (2.6) there) the Arak process ${\cal A}_D$ coincides
 with the Gibbsian modification of the process on $\Gamma_D$ induced by $\Pi_D$
 with the Hamiltonian
 \begin{equation}\label{PHII}
   \Phi(\gamma) := 2 \lgth(\gamma) + \sum_{e \in E(\gamma)} \log \lgth(e)
   - \sum_{x \in V(\gamma)} \log |\sin \phi_x|,
 \end{equation}
 where $E(\gamma)$ and $V(\gamma)$ are, respectively, the edge and vertex sets
 of $\gamma$ while $\phi_x$ stands for the angle between the edges meeting
 in $x$ if $x \in D$ and for the angle between the edge and the tangent to
 $\partial D$ at $x$ if $x \in \partial D.$ This means that
 \begin{equation}\label{GREPR2}
  {\Bbb P}({\cal A}_D \in G) = \frac{{\Bbb E} \sum_{\gamma \in \Gamma_D(\Pi_D) \cap G}
            \exp(-\Phi(\gamma))}{{\Bbb E}\sum_{\gamma \in \Gamma_D(\Pi_D)}
                                      \exp(-\Phi(\gamma))}
 \end{equation}
 for all Borel $G \subseteq \Gamma_D.$
 The third equivalent description of polygonal Markov fields is available in terms
 of equilibrium evolution of one-dimensional particle systems, tracing the polygonal
 realisations of the process in two-dimensional time-space. This description, usually
 referred to as the dynamic representation and introduced already in the original Arak
 work (1982), turned out to be very useful in establishing the essential properties
 of the models. Below, we discuss the dynamic representation for the Arak process,
 see Section 4 in Arak \& Surgailis (1989). We interpret the open convex domain $D$
 as a set of {\it time-space} points $(t,y) \in D,$ with $t$ referred to as the
 {\it time} coordinate and with $y$ standing for the {\it spatial} coordinate of
 a particle at the time $t.$ In this language, a straight line segment in $D$
 stands for a piece of the time-space trajectory of a freely moving particle.
 For a straight line $l$ non-parallel to the time axis and crossing the domain
 $D$ we define in the obvious way its entry point to $D,\;\iin(l,D) \in
 \partial D$ and its exit point $\oout(l,D) \in \partial D.$

 We choose the time-space birth coordinates for the new particles according
 to a homogeneous intensity $\pi$ Poisson point process in $D$ (interior birth
 sites) superposed with a Poisson point process on the boundary (boundary
 birth sites) with the intensity measure
 \begin{equation}\label{KAPPAWPRO}
  \kappa(B) = {\Bbb E}\card\{ l \in \Lambda,\; \iin(l,D) \in B \},\; B \subseteq \partial D.
 \end{equation}
 Each interior birth site emits two particles, moving with initial velocities
 $v'$ and $v''$ chosen according to the joint distribution
 \begin{equation}\label{THETAWPRO}
   \theta(dv',dv'') := \pi^{-1} |v'-v''| (1+{v'}^{2})^{-3/2} (1+{v''}^{2})^{-3/2} dv' dv''.
 \end{equation}
 This can be shown to be equivalent to choosing the directions of the straight lines
 representing the space-time trajectories of the emitted particles according to the
 distribution of the {\it typical angle} between two lines of $\Lambda,$
 see Sections 3 and 4 in Arak \& Surgailis (1989) and the references therein.
 Each boundary birth site $x \in \partial D$ yields one particle
 with initial speed $v$ determined according to the distribution
 $\theta_x(dv)$ identified by requiring that the direction of the line
 entering $D$ at $x$ and representing the time-space trajectory of the
 emitted particle be chosen according to the distribution of a straight
 line $l \in \Lambda$ conditioned on the event $\{ x = \iin(l,D) \}.$

 All the  particles evolve independently in time according to the following rules.
 \begin{description}
  \item{\bf (E1)} Between the critical moments listed below each particle
                  moves freely with constant velocity so that $dy = v dt,$
  \item{\bf (E2)} When a particle touches the boundary $\partial D,$ it dies,
  \item{\bf (E3)} In case of a collision of two particles (equal spatial coordinates $y$
        at some moment $t$ with $(t,y) \in D$), both of them die,
  \item{\bf (E4)} The time evolution of the velocity $v_t$ of an individual particle
        is given by a pure-jump Markov process so that
        $$ {\Bbb P}(v_{t+dt} \in du \;|\; v_t = v) = q(v,du) dt $$
        for the transition kernel
        $$ q(v,du) := |u-v| (1+u^2)^{-3/2} du dt. $$
 \end{description}
 It has been proven (see e.g. Lemma 4.1 in Arak \& Surgailis (1989)) that with
 the above construction of the interacting particle system, the time-space
 trajectories traced by the evolving particles coincide in distribution with
 the Arak process ${\cal A}_D.$  Moreover, a much broader class of consistent
 polygonal Markov fields admit analogous dynamic representations, possibly
 enhanced to allow vertices of higher degrees ($3$ and $4$), see ibidem.
 The question of characterising the class of all polygonal Markov fields
 admitting dynamic representation is far from being trivial and a
 conjectured description of this class has been provided in Arak, Clifford
 \& Surgailis (1993).

 The above dynamic construction of the Arak process makes it very suitable
 for simulation. However, in the present paper we focus our interest on
 the family of processes $\hat{\cal A}_D^{[\alpha,\beta]},\; \alpha, \beta \in {\Bbb R},$
 arising as the Ising-like length- and area-interacting Gibbsian modifications of
 ${\cal A}_D.$ To this end we {\it colour} the original Arak process ${\cal A}_D$
 as follows. Requiring that the polygonal contours of ${\cal A}_D$ stand for
 interfaces between black- and white-coloured regions in $D$ leaves
 us almost surely with two possible ways of colouring $D$ in black
 and white, arising from each other by a simple colour flip. We
 choose one of these colourings at random, with probability $1/2,$
 thus obtaining a coloured version of ${\cal A}_D,$ denoted
 in the sequel by $\hat{\cal A}_D.$ The family of all admissible
 {\it coloured} polygonal configurations in $D,$ carrying information
 not only about the planar contours it consists of, but also about the
 associated colouring, will be denoted by $\hat{\Gamma}_D.$
 With this notation and terminology we define the
 (coloured) processes $\hat{\cal A}_D^{[\alpha,\beta]}$ by
 \begin{equation}\label{GIBBSARAK}
  \frac{d{\cal L}(\hat{\cal A}^{[\alpha,\beta]}_D)}{d{\cal L}(\hat{\cal A}_D)}[\hat{\gamma}]
   := \frac{\exp\left(- {\cal H}_D^{[\alpha,\beta]}(\hat{\gamma}) \right)}
      { {\Bbb E} \exp\left(- {\cal H}_D^{[\alpha,\beta]}(\hat{\cal A}_D) \right)},\;
    \hat{\gamma} \in \hat{\Gamma}_D,
 \end{equation}
  with ${\cal L}(\cdot)$ denoting the law of the argument random object
  and with
 \begin{equation}\label{STAREH}
  {\cal H}_D^{[\alpha,\beta]}(\hat{\gamma}) := \alpha A({\rm black}[\hat{\gamma}])
    + \beta \lgth(\hat{\gamma}),
 \end{equation}
 where ${\rm black}[\hat{\gamma}]$ is the black-coloured region in $D$
 for $\hat{\gamma}$ while $A(\cdot)$ stands for the area measure.
 We also write ${\cal A}_D^{[\alpha,\beta]}$ for the contour ensemble of
 $\hat{\cal A}_D^{[\alpha,\beta]},$ with the colours 'forgotten' and,
 likewise, $\gamma$ for the {\it colourless} version of $\hat{\gamma}
 \in \hat{\Gamma}_D.$ Note that using the symmetry between black and
 white and possibly flipping the colours, whenever convenient we may
 assume without loss of generality that $\alpha \geq 0$ (and we do
 so in the proof of Theorem \ref{ISTNIENIE} below).

 Observe that the modifications of the type (\ref{GIBBSARAK}) fall into
 the general setting considered by Arak \& Surgailis (1989) only for
 $\beta \geq 0,$ see Corollary 4.1 there. However, we find it
 natural to admit also negative $\beta'$s since there is no
 obvious {\it infinite temperature} non-interacting field
 available as the reference object for polygonal Markov fields.
 Consequently, in the sequel we will abuse the language by
 referring to large positive values of $\beta$ as to the
 low temperature region, and to small, possibly negative $\beta'$s
 as to the high temperature regime. For $\beta < 0$ one has to
 check that the partition function
 ${\Bbb E} \exp\left(- {\cal H}_D^{[\alpha,\beta]}(\hat{\cal A}_D) \right)$
 is finite. In Corollary \ref{BETAUJEMNE} we show that this is indeed the
 case and, consequently, the definition (\ref{GIBBSARAK}) is correct
 for all $\beta \in {\Bbb R}.$ Clearly, there are no such problems
 for $\alpha,$ since the overall black or white area is
 deterministically bounded by $A(D).$ It should be emphasised
 though that at present we are able to establish the existence
 of the thermodynamic limit only for $\beta > 0,$
 see Theorem \ref{ISTNIENIE}.

 Models of the type (\ref{GIBBSARAK}) have recently
 found interest in the physical literature, see Nicholls (2001).
 In particular, it has been argued that they exhibit a phase transition
 similar to that of the planar Ising model, with the low temperature
 phase admitting only finite contour nesting (as rigorously shown in Nicholls (2001)),
 and with the high temperature phase conjectured (not yet proven) to exhibit
 infinite contour nesting.

 Below, we shall also consider versions of the above models with
 empty boundary conditions, arising by conditioning the original model
 on the event of there being no vertices on the boundary, so that
 \begin{equation}\label{PUSTYBRZEG}
  {\cal L}\left(\hat{\cal A}^{[\alpha,\beta]}_{D|\emptyset}\right)
  := {\cal L}\left( \hat{\cal A}^{[\alpha,\beta]}_D \right| \left.
     {\cal A}^{[\alpha,\beta]}_D \cap \partial D = \emptyset \right).
 \end{equation}
 In particular,
 $$ \hat{\cal A}_{D|\emptyset} := \hat{\cal A}^{[0,0]}_{D|\emptyset}. $$
 Likewise, we shall consider versions of these models with black
 (or white) boundary conditions given by
 \begin{equation}\label{CZARNYBRZEG}
  {\cal L}\left(\hat{\cal A}^{[\alpha,\beta]}_{D|{\rm black (white)}}\right)
  := {\cal L}\left( \hat{\cal A}^{[\alpha,\beta]}_D \right| \left.
     {\cal A}^{[\alpha,\beta]}_D \cap \partial D = \emptyset,
      \; \partial D \mbox{ is black (white) } \right)
 \end{equation}
 with
 $$ \hat{\cal A}_{D|{\rm black (white)}} := \hat{\cal A}^{[0,0]}_{D|{\rm black (white)}}. $$
 As a direct conclusion from (\ref{GIBBSARAK}) we get
 \begin{equation}\label{GIBBSARAK2}
  \frac{d{\cal L}(\hat{\cal A}^{[\alpha,\beta]}_{D|{\rm bd}})}
       {d{\cal L}(\hat{\cal A}_{D|{\rm bd}})}[\hat{\gamma}]
    = \frac{\exp\left(- {\cal H}^{[\alpha,\beta]}_D(\hat{\gamma})\right)}
      { {\Bbb E} \exp\left(-{\cal H}^{[\alpha,\beta]}_D(\hat{\cal A}_{D|{\rm bd}})\right)},\;
      \hat{\gamma} \in \hat{\Gamma}_D,\; \gamma \cap \partial D = \emptyset
 \end{equation}
 for ${\rm bd} \in \{ \emptyset, {\rm black}, {\rm white} \}. $
 Observe that, unlike the unconditioned finite-volume fields $\hat{\cal A}^{[\alpha,\beta]}_D,\;
 \alpha \neq 0,$ the conditioned fields with monochromatic boundary conditions are well
 defined also for non-convex bounded open $D$ with piecewise smooth boundary. Indeed,
 take any bounded open convex set $D'$ containing $D$ and set
 $\hat{\cal A}^{[\alpha,\beta]}_{D|{\rm bd}},\; {\rm bd} \in \{ {\rm black}, {\rm white} \},$
 to coincide with $\hat{\cal A}^{[\alpha,\beta]}_{D'}$ conditioned on the event that
 no edge hits $\partial D$ and that the colour on $\partial D$
 agrees with that specified by ${\rm bd}.$ The Markov property of polygonal fields
 (see Arak \& Surgailis (1989))  implies that this construction does not depend
 on the choice of $D'.$ Note that this argument does not apply for
 the empty boundary condition ${\rm bd} = \emptyset,$ unless $\alpha=0.$

 The purpose of this paper is to construct for $\alpha, \beta \in {\Bbb R}$
 a family of random dynamics on $\hat{\Gamma}_D$ which leave the distribution
 of $\hat{\cal A}_D^{[\alpha,\beta]}$ invariant. This yields simulating
 algorithms for $\hat{\cal A}_D^{[\alpha,\beta]},$ both of the Metropolis
 type and of perfect type in the spirit of Fern\'andez, Ferrari \& Garcia (1998,2002).
 While the Metropolis algorithm is given for all $\alpha, \beta \in {\Bbb R}$
 and can be readily implemented (which is a subject of our work in progress),
 the perfect scheme is restricted to $\alpha = 0$ and seems more difficult
 to implement, yet its value lies mainly in that it provides
 important theoretical information about the thermodynamic limit behaviour of
 ${\cal A}^{[0,\beta]}$ in the low temperature region  (for large $\beta$)
 and in that it can be used to simulate in finite windows directly from the
 thermodynamic limit. The finite volume dynamics are discussed in the next
 Section \ref{SKOBJ}. In Section \ref{NIESKOBJ} we discuss infinite-volume
 thermodynamic limits of polygonal fields and establish their existence.
 For $\alpha = 0$ and $\beta$ large enough one of our dynamics, constructed
 in Subsection \ref{CBD} below, admits an infinite-volume extension
 and, as mentioned above, it yields a perfect simulation scheme which enables
 us to show in Section \ref{DOSKO} that for ${\cal A}^{[0,\beta]}$ there
 exists exactly one thermodynamic limit without infinite chains, as made
 specific below, and that this limit is isometry invariant
 as well as exponentially $\alpha$-mixing. In particular, it follows that in the class
 of infinite-volume measures without infinite chains there exist exactly
 two extremal infinite-volume Gibbs measures for $\hat{\cal A}^{[0,\beta]},$
 the black-dominated and white-dominated phase, corresponding to the same
 contour distribution. In this context it should be noted that this simple
 picture does not seem to extend to the whole simplex of infinite-volume
 Gibbs measures for ${\cal A}^{[0,\beta]}$: we conjecture the existence
 and sketch, in Section \ref{NIESKOBJ} below, a tentative construction
 of an infinite number of infinite-volume states admitting infinite chains
 and breaking both the translational and rotational symmetry.

  As already mentioned above, the implementation of the algorithms
  described in this paper is a subject of our current work in
  progress. It should be emphasised that an algorithm for
  simulating polygonal Markov fields, very different than ours,
  has already been given in the literature  by Clifford \&
  Nicholls (1994).


\section{Finite volume dynamics}\label{SKOBJ}
 Below we construct two families of random dynamics which leave invariant
 the laws of the Gibbs-modified polygonal random fields
 $\hat{\cal A}_D^{[\alpha,\beta]}$ in a bounded open convex domain
 $D \subseteq {\Bbb R}^2.$ First of these dynamics, leading to a practically
 feasible and easy to implement Metropolis-type simulation algorithm,
 is based on the dynamic representation of the Arak process.
 The second one relies mainly on the point- and line-based representation
 of general polygonal Markov fields and, after some additional work,
 leads to a graphical construction and a perfect algorithm discussed
 in Section \ref{DOSKO}. We postpone the proof of the finiteness of the
 partition function in (\ref{GIBBSARAK}) to Corollary \ref{BETAUJEMNE} below.

\subsection{Disagreement loop birth and death dynamics}\label{DLBD}
 An important concept below will be that of a {\it disagreement loop}, borrowed
 from Schreiber (2004), Section 2.2. This arises from the dynamic construction
 of the Arak process as provided by the evolution rules {\bf (E1-4)} with
 the corresponding birth rules, see (\ref{KAPPAWPRO}) and (\ref{THETAWPRO}).

 Suppose that we observe
 a particular realisation ${\gamma} \in {\Gamma}_D$ of the
 colourless basic Arak process ${\cal A}_D$ and that we modify the
 configuration by adding an extra birth site $x_0$ to the existing
 collection of birth sites for $\gamma,$ while keeping
 the evolution rules {\bf (E1-4)} for all the particles, including the
 the two newly added ones if $x_0 \in D$ and the single newly
 added one if $x_0 \in \partial D.$ Denote the resulting new random
 (colourless) polygonal configuration by $\gamma \oplus x_0.$ A simple
 yet crucial observation is that for $x_0 \in D$ the symmetric difference
 $\gamma \triangle [\gamma \oplus x_0]$ is almost surely a single loop
 (a closed polygonal curve), possibly self-intersecting and possibly
 chopped off by the boundary. Indeed, this is seen as follows.
 The leftmost point of the loop $\gamma \triangle [\gamma \oplus x_0]$
 is of course $x_0.$ Each of the two {\it new} particles $p_1, p_2$
 emitted from $x_0$ move independently, according to ${\bf (E1-4)},$ each
 giving  rise to a {\it disagreement path}. The initial segments of such a
 disagreement path correspond to the movement of a particle, say  $p_1,$
 before its annihilation in the first collision. If this is a collision
 with the boundary, the disagreement path gets chopped off and terminates
 there. If this is a collision with a segment of the original configuration
 $\gamma$ corresponding to a certain {\it old} particle $p_3,$ the
 {\it new} particle $p_1$ dies but the disagreement path continues
 along the part of the trajectory of $p_3$ which is contained in
 $\gamma$ but not in $\gamma \oplus x_0.$ At some further moment $p_3$
 dies itself in $\gamma,$ touching the boundary or killing another
 particle $p_4$ in $\gamma.$ In the second case, however, this collision
 only happens for $\gamma$ and not for $\gamma \oplus x_0$ so the
 particle $p_4$ survives (for some time) in $\gamma \oplus x_0$
 yielding a further connected portion of the disagreement path for $p_1,$
 which is contained in $\gamma \oplus x_0$ but not in $\gamma$ etc.
 A recursive continuation of this construction shows that the disagreement
 path initiated by $p_1$ consists alternately of connected polygonal subpaths
 contained in $[\gamma \oplus x_0] \setminus \gamma$ (call these {\it positive}
 parts) and in $\gamma \setminus [\gamma \oplus x_0]$ (call these {\it negative}
 parts). Note that this disagreement path is self-avoiding and, in
 fact, it can be represented as the graph of some piecewise linear
 function $t \mapsto y(t).$ Clearly, the same applies for the disagreement path
 initiated by $p_2.$ An important observation is that whenever two {\it positive}
 or two {\it negative} segments of the two disagreement paths hit each other,
 both disagreement paths die at this point and the disagreement loop closes
 (as opposed to intersections of segments of distinct signs which do not
 have this effect). Obviously, if the disagreement loop does not close
 in the above way, it gets eventually chopped off by the boundary.
 We shall write $\Delta^{\oplus}[x_0;\gamma] = \gamma \triangle [\gamma \oplus x_0]$
 to denote the (random) disagreement loop constructed above. It remains
 to consider the case $x_0 \in \partial D,$ which is much simpler
 because there is only one particle emitted and so $\Delta^{\oplus}[x_0;\gamma]
 = \gamma \triangle [\gamma \oplus x_0]$ is a single self-avoiding
 polygonal path eventually chopped off by the boundary. We abuse
 the language calling such $\Delta^{\oplus}[x_0;\gamma]$ a (degenerate)
 disagreement loop as well.

 Likewise, a disagreement loop arises if we {\it remove} one birth site $x_0$
 from the collection of birth sites of an admissible polygonal configuration
 $\gamma \in \Gamma_D,$ while keeping the evolution rules for
 all the remaining particles. We write $\gamma \ominus x_0$ for the
 configuration obtained from $\gamma$ by removing $x_0$ from the list
 of the birth sites, while the resulting random disagreement loop is denoted
 by $\Delta^{\ominus}[x_0;\gamma]$ so that $\Delta^{\ominus}[x_0;\gamma]
 = \gamma \triangle [\gamma \ominus x_0].$

 With the above terminology we are in a position to describe a random
 dynamics on the coloured configuration space $\hat{\Gamma}_D,$
 which leaves invariant the law of the
 basic Arak process $\hat{\cal A}_D.$ Particular care is needed, however,
 to distinguish between the notion of time considered in the dynamic
 representation of the Arak process as well as throughout the construction
 of the disagreement loops above, and the notion of time to be introduced
 for the random dynamics on $\hat{\Gamma}_D$ constructed below. To make
 this distinction clear we shall refer to the former as to the {\it
 representation time} (r-time for short) and shall keep for it the notation
 $t,$ while the latter will be called the {\it simulation time} (s-time
 for short) and will be consequently denoted by $s$ in the sequel.

 Consider the following pure jump birth and death type Markovian
 dynamics on $\hat{\Gamma}_D.$
 \begin{description}
  \item{\bf (DL:birth)} With intensity $[\pi dx + \kappa(dx)] ds$
   set $\gamma_{s+ds} := \gamma_s \oplus x$ for $\kappa$ as in
   (\ref{KAPPAWPRO}), then construct $\hat{\gamma}_{s+ds}$
   by randomly choosing, with probability $1/2,$ either of
   the two possible colourings for $\gamma_{s+ds},$
  \item{\bf (DL:death)} For each birth site $x$ in $\gamma_s$
   with intensity $1$ set $\gamma_{s+ds} := \gamma_s
   \ominus x,$ then construct $\hat{\gamma}_{s+ds}$
   by randomly choosing, with probability $1/2,$ either
   of the two possible colourings for $\gamma_{s+ds}.$
 \end{description}
 If none of the above updates occurs we keep $\hat{\gamma}_{s+ds} =
 \hat{\gamma}_s.$ It is convenient to perceive the above dynamics
 in terms of generating random disagreement loops $\lambda$ and
 setting $\gamma_{s+ds} := \gamma_s \triangle \lambda,$ with the
 loops of the type $\Delta^{\oplus}[\cdot,\cdot]$ corresponding
 to the rule {\bf (DL:birth)} and $\Delta^{\ominus}[\cdot,\cdot]$
 to the rule {\bf (DL:death)}.

 As an direct consequence of the dynamic representation
 of the Arak process $\hat{\cal A}_D$ we obtain
 \begin{proposition}\label{AR1}
  The distribution of the Arak process $\hat{\cal A}_D$ is
  the unique invariant law of the dynamics given by {\bf (DL:birth)}
  and {\bf (DL:death)}. The resulting stationary process is reversible.
  Moreover, for any initial distribution of $\hat{\gamma}_0$ the
  laws of the random polygonal fields $\hat{\gamma}_s$ converge in variational
  distance to the law of $\hat{\cal A}_D$ as $s \to \infty.$
 \end{proposition}
 The uniqueness and convergence statements in the above proposition
 require a short justification. They both follow by the observation
 that, in finite volume, regardless of the initial state, the process
 $\hat{\gamma}_s$ spends a non-null fraction of time in the state
 'black' (no contours, the whole domain $D$ coloured black).
 Indeed, this observation allows us to conclude the required uniqueness
 and convergence by a standard coupling argument.

 Below, we  show that the laws of the Gibbs-modified polygonal fields
 $\hat{\cal A}^{[\alpha,\beta]}_D$ arise as the unique invariant
 distributions for appropriate modifications of the reference
 dynamics {\bf (DL:birth), (DL:death)}. The main change is that
 the birth and death updates are no more performed unconditionally,
 they pass an {\it acceptance test} instead and are accepted
 with certain state-dependent probabilities, upon failure of
 the acceptance test the update is discarded. For $a \geq 0,
 b \geq 0$ and $\alpha + a \geq 0, \beta + b \geq 0$ consider
 the following dynamics
 \begin{description}
  \item{${\bf (DL:birth[\alpha,\beta;a,b])}$}
   With intensity $[\pi dx + \kappa(dx)] ds$ do
   \begin{itemize}
    \item put $\delta := \gamma_s \oplus x,$
    \item construct $\hat{\delta}$ by randomly choosing,
      with probability $1/2,$ either of the two possible
      colourings for $\delta,$
    \item accept $\hat{\delta}$ with probability
          $$ \exp\left( - \alpha A\left({\rm black}[\hat{\delta}] \setminus
                             {\rm black}[\hat{\gamma}_s]\right)
               - \beta \lgth( \delta \setminus \gamma_s) \right) $$
          $$ \exp\left( - a A\left({\rm black}[\hat{\delta}] \triangle
                             {\rm black}[\hat{\gamma}_s]\right)
               - b \lgth( \delta \triangle \gamma_s) \right), $$
    \item if accepted, set $\hat{\gamma}_{s+ds} := \hat{\delta},$ otherwise
          keep $\hat{\gamma}_{s+ds} := \hat{\gamma}_s.$
   \end{itemize}
  \item{${\bf (DL:death[\alpha,\beta;a,b])}$}
   For each birth site $x$ in $\gamma_s$ with intensity $1$ do
   \begin{itemize}
    \item put $\delta := \gamma_s \ominus x,$
    \item construct $\hat{\delta}$ by randomly choosing, with probability $1/2,$
      either of the two possible colourings for $\delta,$
    \item accept $\hat{\delta}$ with probability
        $$ \exp\left( - \alpha A\left({\rm black}[\hat{\delta}] \setminus
                             {\rm black}[\hat{\gamma}_s]\right)
               - \beta \lgth( \delta \setminus \gamma_s) \right) $$
          $$ \exp\left( - a A\left({\rm black}[\hat{\delta}] \triangle
                             {\rm black}[\hat{\gamma}_s]\right)
               - b \lgth( \delta \triangle \gamma_s) \right), $$
    \item if accepted, set $\hat{\gamma}_{s+ds} := \hat{\delta},$ otherwise
          keep $\hat{\gamma}_{s+ds} := \hat{\gamma}_s.$
  \end{itemize}
 \end{description}
 In analogy with its original reference form {\bf (DL:birth), (DL:death)},
 the above dynamics should be thought of as generating random disagreement
 loops $\lambda$ and setting $\gamma_{s+ds} := \gamma \triangle \lambda$
 provided $\lambda$ passes the acceptance test. It should be emphasised
 that the random disagreement loops above are generated according to
 the dynamic representation of the original Arak process ${\cal A}_D.$
 The following theorem justifies the above construction.
 \begin{theorem}\label{SYMULACJADL}
  For each $a \geq 0,b \geq 0$ and $\alpha + a \geq 0, \beta + b \geq 0$
  the law of the Gibbs-modified Arak process
  $\hat{\cal A}^{[\alpha,\beta]}_D$ is the unique invariant
  distribution of the dynamics ${\bf (DL:birth[\alpha,\beta;a,b]),
  (DL:death[\alpha,\beta;a,b])}.$ The resulting stationary process is reversible.
  For any initial distribution of $\hat{\gamma}_0$ the laws of the random polygonal
  fields $\hat{\gamma}_s$ converge in variational distance to the law of
  $\hat{\cal A}^{[\alpha,\beta]}_D$ as $s \to \infty.$
 \end{theorem}
  Theorem \ref{SYMULACJADL} can be easily concluded from Proposition \ref{AR1}
  by a straightforward check of the detailed balance conditions. We chose,
  however, to provide below a geometric proof of this result for the
  case $\alpha, \beta \geq 0,$ revealing, in our opinion, the geometric
  intuition underlying the dynamics (a similar proof can be provided
  for $\alpha < 0$ or $\beta < 0$ as well). Note that
  the reason for introducing the additional parameters $a$ and $b$ with the
  possibility that $a > 0, b > 0, \alpha + a > 0$ and $\beta + b > 0$ was
  to gain direct control over the diameter of the region affected by
  a single update, which decays exponentially in the current dynamics.
  The control of the diameter
  of the affected region is a condition {\it sine qua non} for possible
  infinite volume extensions of the ${\bf (DL:...[\alpha,\beta;a,b])}$
  dynamics, which is the subject of our current work in progress.
  Clearly, we could also have chosen another standard set of acceptance
  probabilities conforming to the detailed balance conditions, e.g.
  we could accept a transition $\hat\gamma_s \mapsto
  \hat\gamma_{s+ds} := \hat\delta$
  with probability
  $ \min\left(1,\exp({\cal H}^{[\alpha,\beta]}_D(\hat\gamma_s) -
  {\cal H}^{[\alpha,\beta]}_D(\hat\delta))\right)$ and a direct check of
  the detailed balance conditions, based on Proposition \ref{AR1}, would
  show that the law of $\hat{\cal A}^{[\alpha,\beta]}_D$ is invariant
  with respect to such a dynamics. However, in this dynamics, in
  general we cannot efficiently control the size of the region
  affected in a single update.

  Versions of the disagreement loop birth and death dynamics can be
  easily constructed which leave invariant the distributions of the
  polygonal fields $\hat{\cal A}^{[\alpha,\beta]}_{D|\emptyset}$
  $(\hat{\cal A}^{[\alpha,\beta]}_{D|{\rm black}}, \hat{\cal A}^{[\alpha,\beta]}_{D|{\rm white}})$
  with empty (black,white) boundary conditions respectively.
  To this end, we modify accordingly the dynamics ${\bf (DL:birth[\alpha,\beta;a,b])}$
  and ${\bf (DL:death[\alpha,\beta;a,b])}$ by discarding all the updates
  which make the contour collection $\gamma_s$ hit the boundary,
  and for the monochromatic black or white boundary condition, in addition,
  upon an update we do not pick the colouring by random but we choose
  the unique one compatible with the boundary condition. Denoting the
  so constructed dynamics by ${\bf (DL_{\emptyset}:...[\alpha,\beta;a,b])},$
  ${\bf (DL_{\rm black}:...[\alpha,\beta;a,b])}$ and
  ${\bf (DL_{\rm white}:...[\alpha,\beta;a,b])}$ respectively, we immediately
  conclude the following corollary from Theorem \ref{SYMULACJADL}
 \begin{corollary}\label{SYMULACJADL2}
  For each $a \geq 0 ,b \geq 0$ and $\alpha +a \geq 0, \beta + b \geq 0,$
  the law of the Gibbs-modified Arak process
  $\hat{\cal A}^{[\alpha,\beta]}_{D|\emptyset} (\hat{\cal A}^{[\alpha,\beta]}_{D|{\rm black}},
  \hat{\cal A}^{[\alpha,\beta]}_{D|{\rm white}})$ is the unique invariant distribution
  distribution of the dynamics ${\bf (DL_{\emptyset}:...[\alpha,\beta;a,b])}$
  [${\bf (DL_{\rm black}:...[\alpha,\beta;a,b])}$
   or ${\bf (DL_{\rm white}:...[\alpha,\beta;a,b])})$ respectively].
  The resulting stationary processes are reversible. For any initial distribution of
  $\hat{\gamma}_0$ the laws of the random polygonal fields $\hat{\gamma}_s$
  converge in variational distance to the law of $\hat{\cal A}^{[\alpha,\beta]}_{D|\emptyset}$
  [$\hat{\cal A}^{[\alpha,\beta]}_{D|{\rm black}}, \hat{\cal A}^{[\alpha,\beta]}_{D|{\rm white}}$
   respectively], as $s \to \infty.$
 \end{corollary}

 We believe that a very similar dynamics could be used to simulate length-
 and area-interacting modifications of more general consistent polygonal
 Markov fields admitting the dynamic representation as discussed in
 Arak \& Surgailis (1989,1991) and Clifford, Arak \& Surgailis (1993).
 The only change would be an appropriate redefinition of the operations
 $\Delta^{\oplus}[\cdot;\cdot]$ and $\Delta^{\ominus}[\cdot;\cdot],$
 and the resulting disagreement field would no more be a single
 loop.

\subsection{Contour birth and death dynamics}\label{CBD}
 As already mentioned, unlike the previous one, the dynamics discussed in this subsection
 is constructed in a much narrower setting, restricted to colourless contour
 configurations which do not hit the boundary, and it is meant to leave invariant the
 distributions of ${\cal A}_{D|\emptyset}^{[0,\beta]}.$ Recall from the discussion
 following (\ref{GIBBSARAK2}) that in this setting we can take $D$ to be an arbitrary
 bounded open set in ${\Bbb R}^d,$ with piecewise smooth boundary, we do not need convexity.
 The approach developed
 in this section leads to a simulation algorithm discussed in Section \ref{DOSKO} below,
 which, though perfect, seems to be practically infeasible due to non-constructive description
 of the intensity measure of contour births. However, its value lies in that its infinite
 volume extension provides important theoretical information about the thermodynamic
 limit ${\cal A}^{[0,\beta]},$ yielding in particular the uniqueness of the thermodynamic
 limit for $\beta$ large enough. Observe that the dynamics constructed in this section
 could be in principle also used directly for Metropolis sampling, yet the previous
 disagreement loop dynamics seems much better suited for this particular purpose.

 To proceed, we consider the space ${\cal C}_D$ consisting of all closed polygonal
 contours in $D$ which do not touch the boundary $\partial D.$ For a given point
 configuration $\bar{x} := \{ x_1, ... , x_n \}$ we denote by ${\cal C}_D(\bar{x})$
 the family of those polygonal contours in ${\cal C}_D$ which belong to $\Gamma_D(\bar{x}),$
 i.e. whose vertex set coincides with $\bar{x}.$ We construct the so-called free contour
 measure $\Theta_D$ on ${\cal C}_D$ by putting for $C \subseteq {\cal C}_D$ measurable, say,
 with respect to the Borel $\sigma$-field generated by the Hausdorff distance topology,
 \begin{equation}\label{WOLNEKONTURY}
  \Theta_D(C) := \int_{{\rm Fin}(D)} \sum_{\theta \in C \cap {\cal C}_D(\bar{x})}
   \exp(-\Phi(\theta)) \nu^*(d\bar{x})
 \end{equation}
 with the Hamiltonian $\Phi$ as in (\ref{PHII}), with ${\rm Fin}(D)$ standing
 for the family of finite point configurations in $D$ and where $\nu^*$ is the
 measure on ${\rm Fin}(D)$ given by $d \nu^*(\bar{x}) :=  dx_1 ... dx_n.$
 In order to provide an
 alternative line- rather than point-based expression for $\Theta_D,$ for
 a given finite configuration $(l) := ( l_1,...,l_n )$ of straight lines
 intersecting $D$ denote by ${\cal C}_D(l)$ the family of those polygonal
 contours in ${\cal C}_D$ which belong to $\Gamma_D(l).$ Then we have,
 see e.g. (3.8) in the proof of Theorem 1 in Arak, Clifford \& Surgailis (1993),
 \begin{equation}\label{WOLNEKONTURY2}
   \Theta_D(C) = \int_{{\rm Fin}(L[D])} \sum_{\theta \in C \cap {\cal C}_D(l)}
    \exp(-2\lgth(\theta)) d\mu^*((l))
 \end{equation}
 with ${\rm Fin}(L[D])$ standing for the for the family of finite line
 configurations intersecting $D$ and where $\mu^*$ is the measure
 on ${\rm Fin}(L[D])$ given by $d\mu^*((l_1,...,l_n)) := 
 d\mu(l_1) ... d\mu(l_n)$ with $\mu$ defined in the discussion preceding (\ref{GREPR1}).

 For $\beta \in {\Bbb R}$ we consider the exponential modification $\Theta^{[\beta]}_D$
 of the free measure $\Theta_D,$ given by
 \begin{equation}\label{THETAB}
  \Theta^{[\beta]}_D(d\theta) := \exp(-\beta \lgth(\theta)) \Theta_D(d\theta).
 \end{equation}
 It is easily seen that the total mass $\Theta_D^{[\beta]}({\cal C}_D)$ is
 always finite. Indeed, using (\ref{WOLNEKONTURY2}), taking into account
 that the length of a line segment in $D$ can be at most $\diam(D)$ and
 recalling that, by standard integral geometry,
 $M := \mu(\{ l \;|\; l \cap D \neq \emptyset \}) \leq \lgth(\partial \conv(D))$ we
 conclude that
 $$
  \Theta^{[\beta]}_D({\cal C}_D) \leq \sum_{k=0}^{\infty}
  \frac{M^k \exp(k |\beta| \diam(D))}{k!} \leq  $$
 \begin{equation}\label{SKONCZONOSC}
  \exp[\lgth(\partial \conv(D)) \exp(|\beta|\diam(D)) ] < \infty.
 \end{equation}
 Let ${\cal P}_{\Theta_D^{[\beta]}}$ be the Poisson point process on
 ${\cal C}_D$ with intensity measure $\Theta^{[\beta]}_D.$ It then
 follows directly by (\ref{WOLNEKONTURY}), by the point-based
 representation (\ref{GREPR2}) and by (\ref{PUSTYBRZEG}) that for
 all $\beta \in {\Bbb R}$ for which the partition function
 ${\Bbb E} \exp\left(- {\cal H}_D^{[\alpha,\beta]}(\hat{\cal A}_{D|\emptyset})\right)$
 in (\ref{GIBBSARAK2}) is finite (in fact, we show that this holds for all
 $\beta \in {\Bbb R}$ in Corollary \ref{BETAUJEMNE} below), the polygonal
 field  ${\cal A}^{[0,\beta]}_{D|\emptyset}$ coincides in distribution with
 the union of contours in ${\cal P}_{\Theta_D^{[\beta]}}$
 conditioned on the event that they are disjoint so that
 \begin{equation}\label{WARUNKOWKT}
  {\cal L}\left({\cal A}^{[0,\beta]}_{D|\emptyset}\right) =
  {\cal L}\left( \bigcup_{\theta \in {\cal P}_{\Theta_D^{[\beta]}}} \theta \; \left|
  \; \forall_{\theta, \theta' \in {\cal P}_{\Theta_D^{[\beta]}}} \theta \neq \theta'
  \Rightarrow \theta \cap \theta' = \emptyset \right. \right),
 \end{equation}
 where the conditioning is well defined in view of (\ref{SKONCZONOSC}).
 In particular, taking into account (\ref{GREPR1}) and (\ref{WOLNEKONTURY2}),
 we have for all $\beta$ where (\ref{GIBBSARAK2}) makes sense
 $$ {\Bbb P}\left(\forall_{\theta, \theta' \in {\cal P}_{\Theta_D^{[\beta]}}}
    \theta \neq \theta' \Rightarrow \theta \cap \theta' = \emptyset \right) =
    {\Bbb E}\sum_{\delta \in \Gamma_{D|\emptyset}(\Lambda_D)} \exp(-[2+\beta] \lgth(\delta)) $$
 with $\Gamma_{D|\emptyset}$ standing for the family of admissible polygonal configurations
 in $D$ which do not touch $\partial D.$
 It easily follows that the law of ${\cal A}^{[0,\beta]}_{D|\emptyset}$ is invariant
 and reversible with respect to the following contour birth and death dynamics
 $(\gamma_s)_{s \geq 0}$ on $\Gamma_{D|\emptyset}.$
 \begin{description}
  \item{${\bf (C:birth[\beta])}$} With intensity $\Theta^{[\beta]}_D(d\theta) ds$ do
   \begin{itemize}
    \item Choose a new contour $\theta,$
    \item If $\theta \cap \gamma_s = \emptyset,$ accept $\theta$ and
          set $\gamma_{s+ds} := \gamma_s \cup \theta,$
    \item Otherwise reject $\theta$ and keep $\gamma_{s+ds} := \gamma_s,$
   \end{itemize}
  \item{${\bf (C:death[\beta])}$}
      With intensity $1$ for each contour $\theta \in \gamma_s$
      remove $\theta$ from $\gamma_s$ setting $\gamma_s := \gamma_s \setminus \theta.$
 \end{description}
 It is worth noting that, should we accept all the new-coming contours
 without the disjointness test in the above dynamics, we would get the
 Poisson contour process ${\cal P}_{\Theta_D^{[\beta]}}$ as the
 stationary state.

 Observing that the process $\gamma_s$ constructed above spends a non-null fraction of time in the
 state $\emptyset$ and using a standard coupling argument we are led to
 \begin{theorem}\label{KONTURKI}
  The law of the Gibbs-modified Arak process ${\cal A}^{[0,\beta]}_{D|\emptyset}$
  is the unique invariant distribution of the dynamics ${\bf (C:birth[\beta]),
  (C:death[\beta])}.$ The resulting stationary process is reversible.
  For any initial distribution of $\gamma_0$ the laws of random polygonal
  fields $\gamma_s$ converge in variational distance to
  the law of ${\cal A}^{[0,\beta]}_{D|\emptyset}$ as $s \to \infty.$
 \end{theorem}
  All our results in this section were conditional on the partition function
  in (\ref{GIBBSARAK2}) being finite. We claim here that this holds for all
  $\beta \in {\Bbb R}.$ Indeed, since $\beta = 0$ clearly satisfies this condition as
  corresponding to the basic Arak process ${\cal A}_{D|\emptyset},$
  Theorem  \ref{KONTURKI} can be used for $\beta = 0.$ The dynamics 
  {\bf (C)} above implies that the empty-boundary Arak process
  ${\cal A}_{D|\emptyset}$ is stochastically dominated
  (in the sense of inclusion) by the union of contours in 
  ${\cal P}_{\Theta_D},$ see Corollary \ref{STODOM} below.   
  In particular, by (\ref{SKONCZONOSC}), for all $\alpha,\beta \in {\Bbb R}$
  $$ {\Bbb E} \exp\left(- {\cal H}_D^{[\alpha,\beta]}(\hat{\cal A}_{D|\emptyset}) \right)
     \leq \exp(|\alpha| A(D)) {\Bbb E}\exp\left(|\beta| \sum_{\theta \in {\cal P}_{\Theta_D}}
     \lgth(\theta) \right) = $$
  $$ \exp(|\alpha| A(D)) \exp\left(-\Theta_D({\cal C}_D)\right) \exp\left(\Theta^{[|\beta|]}_D({\cal C}_D)\right) < \infty. $$
  By an appropriate redefinition of $\Theta_D$ admitting edges chopped off by the
  boundary, the same argument can be repeated for ${\cal A}_{D|\emptyset}$
  replaced by ${\cal A}_D.$ Thus, we have proven
  \begin{corollary}\label{BETAUJEMNE}
   For each bounded open domain $D \subseteq {\Bbb R}^2$ both the partition functions
   ${\Bbb E} \exp\left(- {\cal H}_D^{[\alpha,\beta]}(\hat{\cal A}_D) \right)$
   in (\ref{GIBBSARAK})
   and
   ${\Bbb E} \exp\left(- {\cal H}_D^{[\alpha,\beta]}(\hat{\cal A}_{D|\emptyset}) \right)$
   in (\ref{GIBBSARAK2}) are finite for all $\alpha,\beta \in {\Bbb R}.$
  \end{corollary}

\section{Thermodynamic limit}\label{NIESKOBJ}
 The purpose of this section is to define the notion of thermodynamic limit
 for the considered polygonal fields and to establish its existence
 (cf. Surgailis (1991)).

 For a smooth closed simple (non-intersecting) curve $c$ in $D$ by
 the trace of a polygonal configuration $\hat{\gamma}$ on $c,$
 denoted in the sequel by $\hat{\gamma} \wedge c,$ we
 mean the knowledge of
 \begin{itemize}
  \item intersection points and intersection directions of $\hat{\gamma}$
        with $c,$
  \item colouring of points of $c.$
 \end{itemize}
 This concept can be formalised in various compatible ways, yet we keep
 the above informal definition in hope that it does not lead to any
 ambiguities while allowing us to avoid unnecessary technicalities.
 For convenience we assume that no edge of $\hat{\gamma}$ is tangent
 to $c,$ which can be ensured with probability $1$ in view of the
 smoothness of $c.$

 Fix $\alpha, \beta \in {\Bbb R}.$ In view of the Gibbsian
 representations (\ref{GREPR1}), (\ref{GREPR2}) and (\ref{GIBBSARAK})
 we easily check that for each $c$ as above and with $\hat{\theta}$
 standing for a trace on $c$ there exists a stochastic kernel
 $\hat{\cal A}^{[\alpha,\beta]}_{\innt c}(\cdot | \hat{\theta})$
 with the property that
 \begin{equation}\label{WAARUNK}
  {\cal L}_{\innt c}\left(\hat{\cal A}^{[\alpha,\beta]}_D | \hat{\cal A}^{[\alpha,\beta]}_D \wedge
    c = \hat{\theta}\right) =
  {\cal L}_{\innt c}\left(\hat{\cal A}^{[\alpha,\beta]}_{D|{\rm bd}} |
        \hat{\cal A}^{[\alpha,\beta]}_{D|{\rm bd}} \wedge
    c = \hat{\theta}\right) = \hat{\cal A}_{\innt c}^{[\alpha,\beta]}(\cdot|\hat{\theta})
 \end{equation}
 for all bounded open and convex $D \supseteq \innt c$ and for ${\rm bd} \in \{
 {\rm black}, {\rm white} \},$ where ${\cal L}_{\innt c}$ denotes the law of the
 argument random element restricted to $\innt c$ (the interior of $c$).

 Consider the family $\Gamma_{{\Bbb R}^2}$ of whole-plane admissible polygonal
 configurations, determined by {\bf (P1), (P2)} and {\bf (P4)} ({\bf (P3)} is
 meaningless in this context) and by the requirement of local finiteness (any
 bounded set is hit by at most a finite number of edges). Let $\hat{\Gamma}_{{\Bbb R}^2}$
 be the corresponding collection of black-and-white coloured whole-plane
 admissible polygonal configurations. It is natural to define
 the family ${\cal G}(\hat{\cal A}^{[\alpha,\beta]})$ of infinite volume Gibbs
 measures (thermodynamic limits) for $\hat{\cal A}^{[\alpha,\beta]}$ as the
 collection of all probability measures on $\hat\Gamma_{{\Bbb R}^2}$ with the
 accordingly distributed random element $\hat{\cal A}$ satisfying
 \begin{equation}\label{GRANTERM}
  {\cal L}_{\innt c}\left(\hat{\cal A} | \hat{\cal A} \wedge c = \hat{\theta}\right) =
  \hat{\cal A}_{\innt c}^{[\alpha,\beta]}(\cdot|\hat{\theta}).
 \end{equation}
 In addition, we shall consider the family ${\cal G}_{\tau}(\hat{\cal A}^{[\alpha,\beta]})$
 of isometry invariant measures in ${\cal G}(\hat{\cal A}^{[\alpha,\beta]}).$ Using
 an appropriate relative compactness argument much along the same lines as in
 Schreiber (2004) we will readily get the existence of at least one isometry-invariant
 thermodynamic limit for each $\beta > 0.$ 
 \begin{theorem}\label{ISTNIENIE}
  For all $\alpha \in {\Bbb R}$ and $\beta > 0,$ the family
  ${\cal G}_{\tau}(\hat{\cal A}^{[\alpha,\beta]})$ is non-empty.
 \end{theorem}
 Note that for $\alpha = 0$ and $\beta$ large enough this statement
 follows also by Theorem in Surgailis (1991).  

 In the sequel, we will establish certain uniqueness results for the thermodynamic
 limit in the low temperature region within a particular class of infinite-volume
 measures without infinite contours. However, we do conjecture
 that for $\alpha = 0$ outside this class there exists an infinite number of
 extreme infinite-volume phases breaking both the rotational and translational
 symmetries.
 We briefly
 and informally sketch their tentative construction. For the increasing
 sequence of squares $(-n,n)^2,\; n=1,2,...$ we consider a sequence
 of boundary conditions arising by requiring that a large number
 $C(n)$ of edges hit the left-hand side of $(-n,n)^2$ (with the
 intersection points located more or less uniformly over the edge),
 the same number of edges intersect the opposite right-hand side,
 but no edges hit the upper and lower sides. We believe that by choosing
 an appropriate growth rate for $C(n)$ we can assure that the
 resulting sequence of polygonal fields on $(-n,n)^2$ is uniformly
 tight (e.g. in the topology discussed in the proof of Theorem
 \ref{ISTNIENIE}) and the accumulation points of this sequence
 are thermodynamic limits for ${\cal A}^{[0,\beta]}$ with infinite
 number of infinite {\it left-to-right} polygonal chains. Moreover,
 the expected number of such chains hitting a disk of radius $1$
 should exhibit untempered growth to infinity with the distance of
 the centre of the disk from the origin. We conjecture that such
 untempered thermodynamic limits should exist even for $\beta = 0$
 where, in the language of the dynamic time-space construction of
 the basic Arak process, one could, roughly speaking, have an
 infinite-density cloud of particles born at the time $-\infty.$
 Such constructions are possible due to the fact that, under very
 rapid edge density growth with the distance from the origin, one
 can enforce the situation where the influence of the boundary
 conditions on $\partial (-n,n)^2$ competes on equal rights
 or even dominates the stabilising bulk effects within
 $(-n,n)^2.$ Clearly, such phenomena cannot show up in the
 stationary regime, see Schreiber (2004) for a discussion.

\section{Perfect simulation from thermodynamic limit and exponential mixing}\label{DOSKO}
 The purpose of the section is to study the contour birth and death dynamics
 of Subsection \ref{CBD} in context of the perfect infinite-volume simulation scheme
 as developed by Fern\'andez, Ferrari \& Garcia (1998,2002). This approach
 is valid only for sufficiently large $\beta.$ It yields a perfect algorithm
 for simulating thermodynamic limits in finite windows and it allows
 us as well to conclude certain uniqueness and mixing results for
 the thermodynamic limit in low temperature regime.

 To this end, we observe first that for all bounded open sets $D$ with
 piecewise smooth boundary
 the free contour measures $\Theta_D$ as defined in (\ref{WOLNEKONTURY})
 arise as the respective restrictions to ${\cal C}_D$ of the same measure
 $\Theta = \Theta_{{\Bbb R}^2}$ on ${\cal C} := \bigcup_{n=1}^{\infty}
  {\cal C}_{(-n,n)^2},$ in the sequel referred to as the infinite volume
 free contour measure. Indeed, this follows easily by the observation
 that $\Theta_{D_1}$ restricted to ${\cal C}_{D_2}$ coincides with
 $\Theta_{D_2}$ for $D_2 \subseteq D_1.$ In the same way we construct
 the infinite-volume exponentially modified measures $\Theta^{[\beta]} = \Theta_{{\Bbb R}^2}^{[\beta]}.$
 The following  result, which is related to the Lemma in the Appendix
 of Nicholls (2001), will be crucial for our further purposes as
 stating exponential decay of the measure $\Theta^{[\beta]}$ with 
 respect to the contour size. 
 \begin{lemma}\label{NICHOLLS}
  For $\beta \geq 2$ we have
  \begin{equation}\label{FKGREENA}
   \Theta^{[\beta]}(\{ \theta\;|\; dx \in \Ver(\theta),\; \lgth(\theta) > R \})
   \leq 4 \pi \exp(-[\beta-2] R) dx.
  \end{equation}
  Moreover, there exists a constant $\varepsilon > 0$ such that, for $\beta \geq 2,$
  \begin{equation}\label{ZSIGMA}
   \Theta^{[\beta]}(\{ \theta\;|\; {\bf 0} \in \innt \theta,\;
   \lgth(\theta) > R \}) \leq \exp(-[\beta-2+\varepsilon \slash 2] R + o(R)).
  \end{equation}
 \end{lemma}
  We note that, in view of (\ref{WARUNKOWKT}) in Section \ref{CBD}, a standard
  Peierls-type argument can be applied to conclude from Lemma \ref{NICHOLLS} that
  there is no infinite contour nesting for ${\cal A}^{[0,\beta]}$ whenever
  $\beta \geq 2.$

 The approach of Fern\'andez, Ferrari \& Garcia (1998, 2002) specialised
 for our  purposes relies on the following graphical construction, briefly
 sketched below, see ibidem for further details. Choose $\beta \geq 2$ large
 enough, as specified below. Define ${\cal F}({\cal C})$ to be the space of
 countable and locally finite collections of contours from ${\cal C},$ with
 the local finiteness requirement meaning that at most a finite number of
 contours can hit a bounded subset of ${\Bbb R}^2.$ Observe that ${\cal F}({\cal C})
 \subseteq \Gamma_{{\Bbb R}^2}$ (there is no equality since
 ${\cal F}({\cal C})$ contains only bounded closed contours while
 $\Gamma_{{\Bbb R}^2}$ also admits infinite polygonal chains).
  On the s-time-space
 ${\Bbb R} \times {\cal F}({\cal C})$ we construct the stationary
 {\it unconstrained (free)} contour birth and death process $(\varrho_s)_{s \in {\Bbb R}}$
 with the birth intensity measure given by $\Theta^{[\beta]}$ and with the
 death intensity $1.$ Note that {\it unconstrained} or {\it free} means here
 that every new-born contour is accepted regardless of whether it hits the union of
 already existing contours or not, moreover we admit negative time here,
 letting $s$ range through ${\Bbb R}$ rather than ${\Bbb R}_+.$ Observe
 also that we need the birth measure $\Theta^{[\beta]}$ to be finite on
 the sets $\{ \theta \in {\cal C} \;|\; \theta \cap A \neq \emptyset \}$
 for all bounded Borel $A \subseteq {\Bbb R}^2$ in order to have the process
 $(\varrho_s)_{s \in {\Bbb R}}$ well defined on ${\Bbb R} \times {\cal F}({\cal C}).$
 By Lemma \ref{NICHOLLS} this is ensured whenever $\beta \geq 2.$
 It is easily seen that, for each $s \in {\Bbb R},$ $\varrho_s$ coincides in
 distribution with the whole-plane Poisson contour process
 ${\cal P}_{\Theta^{[\beta]}}.$

 To proceed, for the free
 process $(\varrho_s)_{s \in {\Bbb R}}$ we perform the following {\it trimming}
 procedure. We place a directed connection from each s-time-space
 instance of a contour showing up in $(\varrho_s)_{s \in {\Bbb R}}$ and denoted
 by $\theta \times [s_0,s_1),$ with $\theta$ standing for the contour and
 $[s_0,s_1)$ for its lifespan, to all s-time-space contour instances
 $\theta' \times [s'_0,s'_1)$ with $\theta' \cap \theta \neq \emptyset,
 s'_0 \leq s_0$ and $s'_1 > s_0.$ In other words, we connect $\theta \times
 [s_0,s_1)$ to those contour instances which may have affected the
 acceptance status of $\theta \times [s_0,s_1)$ in the {\it constrained}
 contour birth and death dynamics {\bf (C)} as dicussed in Subsection \ref{CBD}.
 These connections yield directed chains of s-time-space contour instances,
 we call them the {\it ancestor chains} in the sequel.
 Following Fern\'andez, Ferrari \& Garcia (2002) the union of
 all ancestor chains stemming from a given contour instance is referred to
 as its {\it clan of ancestors}. Using Lemma \ref{NICHOLLS} combined
 with a general technique of stochastic domination by subcritical multitype
 branching processes as discussed in detail in Fern\'andez, Ferrari \& Garcia
 (1998,2002), for $\beta$ large enough we can ensure that all such
 clans of ancestors are a.s. finite and that a single clan size has
 exponentially decaying tail [i.e. the probability that the clan size
 exceeds $R$ is of order $O(\exp(-cR))$ for some $c > 0$].
 In this case we can uniquely determine the acceptance status
 of all the clan members: contour instances with no ancestors are a.s.
 accepted, which automatically and uniquely determines the acceptance
 status of all the remaining members of the clan by recursive application of
 the inter-contour exclusion rule. Discarding the unaccepted contour instances
 leaves us with an s-time-space representation of a stationary evolution
 $(\gamma_s)_{s \in {\Bbb R}}$ on ${\cal F}({\cal C}) \subseteq
 \Gamma_{{\Bbb R}^2}.$ The graphical construction and the argument in
 Fern\'andez, Ferrari \& Garcia (1998, 2002) specialised to our setting yield
 \begin{theorem}\label{FFG}
  Choose $\beta \geq 2$ large enough so that all the ancestor clans
  in the above graphical construction are a.s. finite and
  a single clan size exhibits exponentially decaying tail.
  Then
  \begin{enumerate}
   \item the ${\cal F}({\cal C})$-valued  process $(\gamma_s)_{s \geq 0}$
    given above is well-defined, stationary and reversible,
  \item the stationary distribution ${\cal L}(\gamma_0)$ on ${\cal F}({\cal C})
    \subseteq \Gamma_{{\Bbb R}^2}$ is isometry invariant and belongs to
    ${\cal G}_{\tau}({\cal A}^{[0,\beta]}),$
  \item the dynamics of $(\gamma_s)_{s \in {\Bbb R}}$ is an infinite-volume
    extension of the contour birth and death dynamics {\bf (C)} as introduced
    in Section \ref{CBD}, i.e. $(\gamma_s)_{s \in {\Bbb R}}$ is a Markov
    process on ${\cal F}({\cal C})$ with the infinitesimal generator
    \begin{equation}\label{GENERATOR}
     [L^{[\beta]} F](\eta) := \int_{\cal C} [F(\eta \cup \{ \theta \}) - F(\eta)]
       d\Theta^{[\beta]}(\theta) + \sum_{\theta \in \eta} [F(\eta \setminus \{ \theta \})
       - F(\eta)]
    \end{equation}
    for $\eta \in {\cal F}({\cal C})$ and bounded $F : {\cal F}({\cal C}) \to
    {\Bbb R}$ such that $F(\eta)$ depends only on $\eta \cap D$ for
    some bounded convex set $D,$
  \item $(\gamma_s)_{s \in {\Bbb R}}$ exhibits exponential s-time-space
    $\alpha$-mixing in that there exists $c > 0$ such that
    $$ \sup_{\stackrel{{\cal E}_1 \in \Im_{B(x,1) \times [s_0,s_1]}}
                      {{\cal E}_2 \in \Im_{B(y,1) \times [s'_0,s'_1]}}}
       |{\Bbb P}({\cal E}_1 \cap {\cal E}_2) - {\Bbb P}({\cal E}_1) {\Bbb P}({\cal E}_2)|
       \leq e^{-c [\dist(x,y) + \dist([s_0,s_1],[s'_0,s'_1])]} $$
    whenever $\dist(x,y)$ is sufficiently large, with $\Im_{B(x,1) \times [s_0,s_1]}$ standing for the
    $\sigma$-field generated by the restriction of $(\gamma_s)_{s \in {\Bbb R}}$ to the
    s-space-time region $B(x,1) \times [s_0,s_1],$ where $B(x,1)$ is the
    disk of radius $1$ centred at $x \in {\Bbb R}^2,$
  \item consequently, the stationary distribution ${\cal L}(\gamma_0)$ exhibits
    exponential spatial $\alpha$-mixing.
 \end{enumerate}
\end{theorem}
 It is worth noting that even if $\beta$ is not large enough to ensure a.s.
 finiteness of ancestor clans, a weaker version of the above graphical
 construction can be provided as soon as the birth intensity measure $\Theta^{[\beta]}$
 is finite on $\{ \theta \in {\cal C} \;|\; \theta \cap A \neq \emptyset \}$
 for all bounded $A \subseteq {\Bbb R}^2,$ which is the case whenever $\beta \geq 2$
 by Lemma \ref{NICHOLLS}. To this end we restrict the s-time
 to ${\Bbb R}_+$ and choose an initial condition, which is an ${\cal F}({\cal C})$-valued
 random element independent of the free birth and death process of the graphical
 construction. The birth and death process here is also restricted to positive
 times in that there are no contours born or alive before the s-time $0,$ in
 other words the birth and death process starts with the initial state $\emptyset$
 at s-time $0,$ consequently it is no more stationary. In this context the
 local finiteness of $\Theta^{[\beta]}$ allows us to conclude that for each
 contour instance $\theta \times [s_0,s_1),\; s_0,s_1 > 0,$ the expected
 cardinality of its ancestor clan extending down to the s-time $0$ is finite,
 consequently the clan is a.s. finite (note that it could extend to an infinite
 clan through negative s-times in the original graphical construction). Thus,
 with the initial state given, the acceptance status of each contour instance
 is uniquely determined by the inter-contour exclusion rule. This leads us to
 \begin{corollary}\label{DODATEK}
  With $\beta \geq 2,$ for each ${\cal F}({\cal C})$-valued initial condition
  $\gamma_0$ there exists a Markov process $(\gamma_s)_{s \geq 0}$ on
  ${\cal F}({\cal C})$ with infinitesimal generator given by (\ref{GENERATOR}).
 \end{corollary}
 In the remaining part of the present section we will not use Corollary
 \ref{DODATEK} and, unless otherwise stated, we shall assume that
 $\beta$ stays within the region of validity of the original graphical
 construction preceding Theorem \ref{FFG}. We denote by $\mu^{[\beta]}$ the
 infinite-volume stationary distribution ${\cal L}(\gamma_0)$ arising in this
 graphical construction. Observe that the fact that $\mu^{[\beta]}$ is
 concentrated on ${\cal F}({\cal C})$ means that it contains no
 infinite polygonal chains \--- all the contours are bounded and closed. Below, we show
 that, with the assumptions of Theorem \ref{FFG}, $\mu^{[\beta]}$ is in fact the unique
 element of ${\cal G}({\cal A}^{[0,\beta]})$ concentrated on ${\cal F}({\cal C}),$ although
 we conjecture that ${\cal G}({\cal A}^{[0,\beta]})$ is non-empty as argued in Section
 \ref{NIESKOBJ}. To proceed with our argument we consider finite-volume versions
 of the above graphical construction, with the infinite-volume
 birth intensity measure $\Theta^{[\beta]}$ replaced by its finite
 volume restrictions $\Theta^{[\beta]}_D$ for bounded and open $D$
 with piecewise smooth boundary.
 Clearly, the graphical construction yields
 then a version of the finite-volume contour birth and death dynamics ${\bf (C)}$ as
 discussed in Subsection \ref{CBD}. For each $D$ denote the resulting finite-volume
 stationary process on ${\cal F}({\cal C}_D)$ by $(\gamma_s^D)_{s \in {\Bbb R}}.$
 Write also $(\varrho_s^D)$ for the corresponding {\it free} contour birth
 and death process. Note that this finite-volume construction is valid for all
 $\beta \in {\Bbb R},$ even though in this section it is only used for
 $\beta$ as in Theorem \ref{FFG}. In view of Theorem
 \ref{KONTURKI} we see that $\gamma_s^D$ coincides in distribution
 with ${\cal A}^{[0,\beta]}_{D|\emptyset}$ for all $s \in {\Bbb R}.$
 Moreover, it is easily seen that $\varrho^D_s$ coincides in distribution
 with ${\cal P}_{\Theta_D^{[\beta]}}$ for all $s \in {\Bbb R}.$ 
 From the construction, Lemma \ref{NICHOLLS} and the general theory developed in
 Fern\'andez, Ferrari \& Garcia (1998, 2002) it follows that
 \begin{proposition}\label{ZEBRANIE}
  With $\beta$ as in Theorem \ref{FFG} the finite-volume graphical constructions
  for different $D \subseteq {\Bbb R}^2$ and the infinite-volume graphical
  construction can be coupled on a common probability space so that there
  exists $c > 0$ with
  $$ {\Bbb P}\left(\gamma^{D_1}_s \cap B(x,1) \neq \gamma^{D_2}_s \cap B(x,1)\right)
     \leq \exp(-c \min(\dist(x,\partial D_1),\dist(x,\partial D_2))) $$
  for bounded $D_1, D_2 \subseteq {\Bbb R}^2,$ for $x$ sufficiently far from
  $\partial D_1$ and $\partial D_2$ and for all $s \in {\Bbb R}.$ Moreover,
  $$ {\Bbb P}\left(\gamma^D_s \cap B(x,1) \neq \gamma_s \cap B(x,1)\right)
     \leq \exp(-c \dist(x,\partial D)) $$
  for bounded $D \subseteq {\Bbb R}^2,$ for $x$ far enough
  from $\partial D$ and for all $s \in {\Bbb R}.$
 \end{proposition}
 Taking into account that, by the construction and by the results of Section \ref{CBD},
 $\gamma^D_s$ coincides in distribution with ${\cal A}^{[0,\beta]}_{D|\emptyset},$ and
 that for each contour collection in ${\cal F}({\cal C})$ every bounded region can be
 surrounded by a smooth curve which does not hit any of the contours, we can use the
 Markov property of the considered polygonal fields combined with Proposition \ref{ZEBRANIE}
 to conclude the claimed property
 \begin{corollary}
  For $\beta$ as in Theorem \ref{FFG} the measure $\mu^{[\beta]}$ is the only element
  of ${\cal G}({\cal A}^{[0,\beta]})$ concentrated on ${\cal F}({\cal C}).$
 \end{corollary}
 For $\beta$ as in Theorem \ref{FFG}, using Lemma \ref{NICHOLLS} we easily conclude
 that the number of contours in $\varrho_0$ surrounding a given point
 is a.s. finite. Consequently, the number of contours surrounding
 a given point in $\gamma_0$ is a.s. finite as well, whence there is no infitite contour
 nesting. Thus, we observe a unique infinite connected region surrounding finitely
 nested contour collections. Colouring this region black or white gives rise to
 two distinct phases, respectively black- and white-dominated. There are no other
 extreme phases without infinite chains in the coloured model, because their
 corresponding colourless contour ensembles have to coincide with $\mu^{[\beta]}.$

 The last important conclusion of the graphical construction, based on the
 above-made observations that almost surely $\gamma_s \subseteq \varrho_s,
 \gamma_s^D \subseteq \varrho_s^D$ and that $\gamma_s \sim^d {\cal A}^{[0,\beta]},
 \gamma_s^D $ $\sim^d {\cal A}^{[0,\beta]}_{D|\emptyset}, \varrho_s \sim^d
 {\cal P}_{\Theta^{[\beta]}}$
 and $\varrho_s^D \sim^d {\cal P}_{\Theta_D^{[\beta]}},$ is the following stochastic
 domination statement
 \begin{corollary}\label{STODOM}
  The Poisson contour process ${\cal P}_{\Theta^{[\beta]}}$ stochastically dominates
  (in the sense of inclusion of contour collections) the polygonal field
  ${\cal A}^{[0,\beta]}.$ Likewise,
  for each bounded $D$ with piecewise smooth boundary, the Poisson process 
  ${\cal P}_{\Theta^{[\beta]}_D}$ stochastically dominates the finite-volume
  polygonal field ${\cal A}^{[0,\beta]}_D.$
 \end{corollary} 


\section{Proofs}

\subsection{Proof of Theorem \ref{SYMULACJADL}}
 In order to provide a geometrical and intuitive proof of the
 theorem we construct an auxiliary model. For $r > 0$ define
 $\hat{\cal A}^{[\alpha,\beta;r]}_D$ to be the Gibbsian
 modification of $\hat{\cal A}_D$ with the Hamiltonian
 \begin{equation}\label{NOWEH}
  {\cal H}^{[\alpha,\beta;r]}_D(\hat{\gamma}) :=
  r^{-1} \beta A\left( [\gamma +_M B(r)] \cap D \right) +
  \alpha A(({\rm black}[\hat{\gamma}] +_M B(r)) \cap D),
 \end{equation}
 with $+_M$ standing for the usual Minkowski addition
 and with $B(r)$ denoting the radius $r$ disk in
 ${\Bbb R}^2,$ centred in $0.$ It is easily seen that,
 for each $\hat{\gamma} \in \hat{\Gamma}_D,$
 \begin{equation}\label{DOZERA}
  \lim_{r\to 0} {\cal H}^{[\alpha,\beta;r]}_D(\hat{\gamma}) =
   {\cal H}^{[\alpha,\beta]}_D(\hat{\gamma})
 \end{equation}
 so that ${\cal H}^{[\alpha,\beta;r]}_D$ is an approximation
 of ${\cal H}^{[\alpha,\beta]}_D$ for small $r.$ Take
 $\Pi^{[\alpha+a]}, \Pi^{[r^{-1}(\beta+b)]},$ $ \Pi^{[a]}$ and $\Pi^{[r^{-1} b]}$
 to be independent homogeneous Poisson point processes on $D,$ jointly independent
 of $\hat{\cal A}_D,$ with respective intensities $\alpha+a,r^{-1}(\beta+b), a$
 and $r^{-1} b.$ We claim that $\hat{\cal A}^{[\alpha,\beta;r]}_D$ coincides in
 distribution with $\hat{\cal A}_D$ conditioned jointly with
 $\Pi^{[\alpha+a]}, \Pi^{[r^{-1}(\beta+b)]},$ $ \Pi^{[a]}$ and $\Pi^{[r^{-1} b]}$
 on the event ${\cal E}[\alpha,\beta;a,b;r]$ that the following conditions are
 simultaneously satisfied
 \begin{itemize}
  \item $\Pi^{[r^{-1}(\beta+b)]} \cap [\gamma +_M B(r)] = \emptyset,$
  \item $\Pi^{[r^{-1} b]} \subseteq [\gamma +_M B(r)],$
  \item $\Pi^{[\alpha+a]} \cap [{\rm black}(\hat{\gamma}) +_M B(r)] = \emptyset,$
  \item $\Pi^{[a]} \subseteq [{\rm black}(\hat{\gamma}) +_M B(r)] = \emptyset,$
 \end{itemize}
 so that
 \begin{equation}\label{WARUNKOW}
  {\cal L}(\hat{\cal A}^{[\alpha,\beta;r]}_D) = {\cal L}\left( \hat{\cal A}_D | {\cal E}[\alpha,\beta;a,b;r]
          \right).
 \end{equation}
 Indeed, for a given $\hat{\gamma} \in \hat{\Gamma}_D$ the probability of
 the event ${\cal E}[\alpha,\beta;a,b;r]$ is
 $$ {\Bbb P}({\cal E}[\alpha,\beta;a,b;r]|\hat{\gamma}) =
    \exp\left(- r^{-1} [\beta+b] A([\gamma +_M B(r)] \cap D)\right) $$
 $$ \exp\left(- r^{-1} b [A(D) - A([\gamma +_M B(r)] \cap D)]\right) $$
 $$ \exp\left(- [\alpha+a] A([{\rm black}(\hat{\gamma}) +_M B(r)] \cap D)\right) $$
 $$ \exp\left(- a [A(D) - A([{\rm black}(\hat{\gamma}) +_M B(r)] \cap D)]\right) = $$
 $$ \exp\left( - {\cal H}^{[\alpha,\beta;r]}_D(\hat{\gamma}) \right) \exp\left(-[a+r^{-1}b] A(D) \right), $$
 which yields (\ref{WARUNKOW}) by the definition of $\hat{\cal A}^{[\alpha,\beta;r]}_D.$

 To proceed, we construct an auxiliary Markovian dynamics which leaves invariant the joint distribution
 of of $\hat{\cal A}_D, \Pi^{[r^{-1} (\beta + b)]}, \Pi^{[\alpha + a]}, \Pi^{[r^{-1}b]}$
 and $\Pi^{[a]},$ and makes the resulting stationary process reversible. To this end, set
 $$ \hat{\gamma}_0 := \hat{\cal A}_D, \; \pi^{\alpha}_0 := \Pi^{[\alpha+a]}, \;
    \pi^{\beta}_0 := \Pi^{[\beta+b]},\; \pi^{a}_0 := \Pi^{[a]}, \; \pi^b_0 := \Pi^{[b]} $$
 and let the quintuple $(\hat{\gamma}_s,\pi^{\alpha}_s,\pi^{\beta}_s,\pi^{a}_s,\pi^b_s)_{s \geq 0}$
 evolve according to the following rules, applied independently for each component,
 \begin{description}
  \item{\bf (Aux1)} $\hat{\gamma}_s$ evolves according to ${\bf (DL:birth)}$ and ${\bf (DL:death)},$
  \item{\bf (Aux2)} $\pi^{\alpha}_s, \pi^{\beta}_s, \pi^a_s$ and $\pi^b_s$ evolve according to a birth
        and death process with death intensity $1$ and with birth intensities $\alpha+a,$
        $r^{-1}(\beta+b),$ $a$ and $r^{-1}b$ respectively.
 \end{description}
 The above invariance and reversibility statements follow as direct consequences of
 Proposition \ref{AR1}. Thus, we conclude that the joint distribution of
 $(\hat{\cal A}_D,\Pi^{[r^{-1} (\beta + b)]},$ $\Pi^{[\alpha+a]},$ $\Pi^{[r^{-1}b]},\Pi^{[a]})$
 conditioned on the event ${\cal E}[\alpha,\beta;a,b;r],$ is invariant and reversible
 with respect to the following Markovian dynamics, arising from {\bf (Aux1)} and {\bf (Aux2)} by
 adding an appropriate acceptance test to be passed only by admissible updates:
 \begin{description}
  \item{\bf (B1)} Choose an update $(\hat{\delta},\theta^{\alpha},\theta^{\beta},
   \theta^a,\theta^b)$ for $(\hat{\gamma}_{s+ds},\pi^{\alpha}_{s+ds},\pi^{\beta}_{s+ds},
   \pi^a_{s+ds},\pi^b_{s+ds})$ according to the rules ${\bf (Aux1), (Aux2)}.$
  \item{\bf (B2)} Accept the update, setting
    $$  (\hat{\gamma}_{s+ds},\pi^{\alpha}_{s+ds},\pi^{\beta}_{s+ds},\pi^a_{s+ds},\pi^b_{s+ds})
     := (\hat{\delta},\theta^{\alpha},\theta^{\beta},\theta^a,\theta^b), $$
    provided the following conditions are satisfied
    \begin{itemize}
     \item $\theta^{\beta} \cap [\delta +_M B(r)] = \emptyset,$
     \item $\theta^b \subseteq [\delta +_M B(r)],$
     \item $\theta^{\alpha} \cap [{\rm black}(\hat{\delta}) +_M B(r)] = \emptyset,$
     \item $\theta^a \subseteq [{\rm black}(\hat{\delta}) +_M B(r)],$
    \end{itemize}
  \item{\bf (B3)} Otherwise discard the update, keeping
   $$ (\hat{\gamma}_{s+ds},\pi^{\alpha}_{s+ds},\pi^{\beta}_{s+ds},\pi^a_{s+ds},\pi^b_{s+ds}) :=
      (\hat{\gamma}_{s},\pi^{\alpha}_{s},\pi^{\beta}_{s},\pi^a_{s},\pi^b_{s}). $$
 \end{description}
 Consequently, in view of (\ref{WARUNKOW}), the first component $\hat{\gamma}_s$
 under the above stationary dynamics ${\bf (B1-3)},$ with the initial distribution at $s=0$ given by the joint law of
 $(\hat{\cal A}_D,\Pi^{[r^{-1} (\beta + b)]},$ $\Pi^{[\alpha+a]},$ $ \Pi^{[r^{-1}b]},$ $\Pi^{[a]})$
 conditioned on the event ${\cal E}[\alpha,\beta;a,b;r],$ coincides in distribution with
 $\hat{\cal A}_D^{[\alpha,\beta;r]}$ for all $s \in {\Bbb R}_+.$ Moreover, the conditional
  distributions of the remaining
 components given $\hat{\gamma}_s$ are also readily determined. Indeed, $\pi^{\alpha}_s$
 is just a homogeneous Poisson point process on $D \setminus [{\rm black}(\hat{\gamma}_s) +_M B(r)]$
 with intensity $\alpha+a,$ $\pi^{\beta}_s$ is an intensity $r^{-1}(\beta+b)$ homogeneous Poisson
 point process on $D \setminus [\gamma +_M B(r)],$ $\pi^a_s$ is a homogeneous
 Poisson point process on ${\rm black}(\hat{\gamma}_s) +_M B(r)$ of intensity $a$ while $\pi^b_s$
 is a homogeneous Poisson point process on $\gamma_s +_M B(r)$ with intensity $r^{-1} b.$ All four
 components $\pi^{\alpha}_s, \pi^{\beta}_s, \pi^a_s, \pi^b_s$ are jointly independent
 given $\hat{\gamma}_s.$ Consequently, we observe that if we integrate out the Poisson
 components $\pi^{\alpha}, \pi^{\beta}, \pi^a$ and $\pi^b,$ the polygonal field component
 $\hat{\gamma}_s$ turns out to evolve according to the following dynamics (see Subsection \ref{DLBD}
  for the notation):
 \begin{description}
  \item{${\bf (DL:birth[\alpha,\beta;a,b;r])}$}
   With intensity $[\pi dx + \kappa(dx)] ds$ do
   \begin{itemize}
    \item put $\delta := \gamma_s \oplus x,$
    \item construct $\hat{\delta}$ by randomly choosing,
      with probability $1/2,$ either of the two possible
      colourings for $\delta,$
    \item accept $\hat{\delta}$ with probability
          $$ \exp\left( - [\alpha+a] A\left([{\rm black}(\hat{\delta})+_M B(r)] \setminus
                                            [{\rm black}(\hat{\gamma}_s) +_M B(r)] \right)\right) $$
          $$ \exp\left( - r^{-1}[\beta+b] A([\delta +_M B(r)] \setminus [\gamma_s +_M B(r)]) \right) $$
          $$ \exp\left( - a A\left([{\rm black}(\hat{\gamma}_s) +_M B(r)] \setminus
                                   [{\rm black}(\hat{\delta}) +_M B(r)]\right)\right) $$
          $$ \exp\left( - r^{-1} b A( [\gamma_s +_M B(r)] \setminus [\delta +_M B(r)]) \right) = $$
          $$ \exp\left( - \alpha A\left([{\rm black}(\hat{\delta})+_M B(r)] \setminus
                             [{\rm black}(\hat{\gamma}_s) +_M B(r)] \right)\right) $$
          $$ \exp\left( - \beta r^{-1} A( [\delta +_M B(r)] \setminus [\gamma_s +_M B(r)]) \right) $$
          $$ \exp\left( - a A\left([{\rm black}(\hat{\delta}) +_M B(r)] \triangle
                             [{\rm black}(\hat{\gamma}_s) +_M B(r)]\right)\right) $$
          $$ \exp\left( - b r^{-1} A( [\delta +_M B(r)] \triangle [\gamma_s +_M B(r)]) \right), $$
    \item if accepted, set $\hat{\gamma}_{s+ds} := \hat{\delta},$ otherwise
          keep $\hat{\gamma}_{s+ds} := \hat{\gamma}_s.$
   \end{itemize}
  \item{${\bf (DL:death[\alpha,\beta;a,b;r])}$}
   For each birth site $x$ in $\gamma_s$ with intensity $1$ do
   \begin{itemize}
    \item put $\delta := \gamma_s \ominus x,$
    \item construct $\hat{\delta}$ by randomly choosing, with probability $1/2,$
      either of the two possible colourings for $\delta,$
    \item accept $\hat{\delta}$ with probability
         $$ \exp\left( - \alpha A\left([{\rm black}(\hat{\delta})+_M B(r)] \setminus
                             [{\rm black}(\hat{\gamma}_s) +_M B(r)] \right)\right) $$
         $$ \exp\left( - \beta r^{-1} A( [\delta +_M B(r)] \setminus [\gamma_s +_M B(r)]) \right) $$
         $$ \exp\left( - a A\left([{\rm black}(\hat{\delta}) +_M B(r)] \triangle
                             [{\rm black}(\hat{\gamma}_s) +_M B(r)]\right)\right) $$
         $$ \exp\left( - b r^{-1} A( [\delta +_M B(r)] \triangle [\gamma_s +_M B(r)]) \right), $$
    \item if accepted, set $\hat{\gamma}_{s+ds} := \hat{\delta},$ otherwise
          keep $\hat{\gamma}_{s+ds} := \hat{\gamma}_s.$
  \end{itemize}
 \end{description}
 Thus, the distribution of $\hat{\cal A}^{[\alpha,\beta;r]}_D$ is invariant
 and reversible with respect to the above dynamics. Moreover, it is easily seen that
 the acceptance probabilities in the rules ${\bf (DL:birth[\alpha,\beta;a,b;r])}$ and
 ${\bf (DL:death[\alpha,\beta;a,b;r])}$ converge to these in ${\bf (DL:birth[\alpha,\beta;a,b])}$
 and ${\bf (DL:death[\alpha,\beta;a,b])}$ as $r \to 0.$ Taking into account (\ref{DOZERA})
 and letting $r \to 0$ we get by a standard continuity argument that $\hat{\cal A}^{[\alpha,\beta]}_D$
 is invariant and reversible with respect to the dynamics ${\bf (DL:birth[\alpha,\beta;a,b])}$
 and ${\bf (DL:death[\alpha,\beta;a,b])}.$

 To complete the proof of Theorem \ref{SYMULACJADL} it suffices now to establish the
 remaining uniqueness and convergence statements. These follow, however, along the
 same lines as in Proposition \ref{AR1}, by the observation that, in finite volume,
 regardless of the initial state, the process $\hat{\gamma}_s$ spends a non-null
 fraction of time in the state 'black' (no contours, the whole domain $D$ coloured black)
 and by a standard application of coupling argument. The proof is complete. $\Box$


\subsection{Proof of Theorem \ref{ISTNIENIE}}
 Following the ideas of Schreiber (2004) it is convenient to consider the family
 $\Gamma_{{\Bbb R}^2}$ of admissible configurations in the plane embedded into
 the space $G_{{\Bbb R}^2}$ of locally finite non-negative Borel measures on
 ${\Bbb R}^2,$ by identifying a configuration $\hat{\gamma} \in \Gamma_{{\Bbb R}^2}$
 with the measure
 \begin{equation}\label{MMM}
  M_{\hat{\gamma}}(U) := \lgth(\gamma \cap U) + A({\rm black}(\hat{\gamma}) \cap U)
  + N(\gamma \cap U)
 \end{equation}
 for Borel $U \subseteq {\Bbb R}^2,$ with $N(\gamma \cap U)$ standing for the
 number of vertices of $\gamma$ falling into $U.$ Endow the space $G_{{\Bbb R}^2}$
 with the vague topology defined as the weakest one to make continuous the
 mappings $\mu \mapsto \int f d\mu$ for all continuous $f$ with bounded support.
 Observe that in general $\Gamma_D \not\subseteq \Gamma_{{\Bbb R}^2}$
 for $G \subset {\Bbb R}^2$ due to the presence of edges chopped off by
 the boundary. Therefore, in order to have our embedding defined also
 for finite-volume configurations, we agree to put $M_{\hat{\gamma}}(D^c) := 0$
 for all $\hat{\gamma} \in \Gamma_D.$ Note that only the internal vertices of
 finite-volume configurations are counted in $N(\cdot).$     
      
 Consider the sequence $((-n,n)^2)_{n=1}^{\infty}$ of growing open
 squares in ${\Bbb R}^2.$ By the properties of the basic Arak process,
 see Section 4 in Arak \& Surgailis (1989) and Section 2.1 in
 Schreiber (2004), it immediately follows that there exists a
 finite constant $C$ with 
 \begin{equation}\label{OG12}
  {\Bbb E} M_{\hat{\cal A}^{[0,0]}_{(-n,n)^2}}((-n,n)^2) \leq C A((-n,n)^2)
 \end{equation}
 for all $n \geq 1.$ We will show that the above conclusion can be extended for arbitrary
 $\alpha \in {\Bbb R}$ and $\beta > 0$ in that there exists $C^{[\alpha,\beta]} < \infty$
 with 
 \begin{equation}\label{OGG12}
  {\Bbb E} M_{\hat{\cal A}^{[\alpha,\beta]}_{(-n,n)^2}}((-n,n)^2) \leq C^{[\alpha,\beta]} A((-n,n)^2).
 \end{equation}
 Below, we assume without loss of generality that $\alpha \geq 0,$ which can be
 done in view of the colour-flip symmetry. Observe first that, in view of (\ref{GIBBSARAK}), 
 $$ \frac{\partial}{\partial h} {\Bbb E} {\cal H}^{[\alpha,\beta]}_{(-n,n)^2}
    (\hat{\cal A}^{[h\alpha,h\beta]})
   = - \Var\left({\cal H}^{[\alpha,\beta]}_{(-n,n)^2}
      (\hat{\cal A}^{[h\alpha,h\beta]}_{(-n,n)^2})\right) < 0 $$
 with ${\cal H}^{[\alpha,\beta]}_{(-n,n)^2}$ as in (\ref{STAREH}). Consequently,
 taking into account that the area term in the Hamiltonian
 ${\cal H}^{[\alpha,\beta]}_{(-n,n)^2}$ is
 bounded by $\alpha A((-n,n)^2)$ and that the Hamiltonian is 
 always positive, we conclude by (\ref{OG12}) that also the expectation
 of the edge length
 term in the Hamiltonian admits an area-order upper bound. It remains to show that this
 is also the case for the number of vertices \--- we sketch the argument omitting standard
 technical details. To this end, we take advantage of the dynamic representation (as
 discussed in the introduction of this paper and in Section 4 of Arak \& Surgailis
 (1989)) to conclude that for the basic Arak process ${\cal A}_{(-n,n)^2}$ the number of internal
 left-extreme vertices (with the corresponding sharp angle lying to the right of the
 vertex) is $\Po(\pi A((-n,n)^2)),$ where $\Po(\tau)$ stands for Poisson-distributed
 random variable with mean $\tau.$
 The same  applies for the number of
 internal right-extreme, upper-extreme and lower-extreme vertices (recall that we do not
 count the boundary vertices here). Consequently, the overall number of internal
 vertices $N({\cal A}_{(-n,n)^2})$ is stochastically bounded by $4\Po(4\pi n^2)$
 and has its mean of area order, not greater than $16 \pi n^2.$ In view of the
 representation (\ref{GIBBSARAK}) and taking into account that the Hamiltonian
 ${\cal H}_{(-n,n)^2}^{[\alpha,\beta]}$ is always positive since $\alpha \geq 0,$
 we conclude that, for all $K > 0,$ 
 \begin{equation}\label{POI}
   {\Bbb P}\left(N(\hat{\cal A}^{[\alpha,\beta]}_{(-n,n)^2}) > 4 K\right) \leq 
   \frac{{\Bbb P}(\Po(4\pi n^2) > K)}{{\Bbb E}
   \exp(-{\cal H}^{[\alpha,\beta]}_{(-n,n)^2}(\hat{\cal A}^{[\alpha,\beta]}_{(-n,n)^2}))}.
 \end{equation}    
 Recall that Poisson distributions exhibit superexponentially decaying tails 
 $$ {\Bbb P}(\Po(4 \pi n^2) > K) \leq \exp\left(-\frac{K}{4} \log\left(\frac{K}{8 \pi n^2}\right) \right),\;
     K \geq 64 \pi n^2,
 $$ 
 see Shorack \& Wellner (1986), p. 485. Moreover, the negative logarithm of the denominator
 in (\ref{POI}) exhibits at most area-order growth, which is due to the easily 
 verified finiteness of the free energy density for ${\cal H}^{[\alpha,\beta]}$
 $$ \liminf_{n\to\infty}
     \frac{1}{(2n)^2} \log {\Bbb E}
     \exp\left(-{\cal H}^{[\alpha,\beta]}_{(-n,n)^2}(\hat{\cal A}_{(-n,n)^2})\right) > - \infty. $$
 Consequently, the required area-order bound for
 ${\Bbb E} N(\hat{\cal A}^{[\alpha,\beta]}_{(-n,n)^2})$ follows now from (\ref{POI}) by
 a direct calculation. This completes the verification of (\ref{OGG12}). 

 To proceed with the proof of the theorem, consider the sequence
 $(M_{n}^{[\alpha,\beta]})_{n=1}^{\infty}$ of
 $G_{{\Bbb R}^2}$-valued random elements with laws given by
 $$
   {\cal L}(M_{n}^{[\alpha,\beta]}) :=
     \frac{1}{4\pi(2n)^2} \int_{[0,2\pi)} \int_{(-n,n)^2}
     {\cal L}\left( [\tau_x \circ R_{\phi}] M_{\hat{\cal A}_{(-n,n)^2}^{[\alpha,\beta]}}\right) dx
     d \phi + $$
 \begin{equation}\label{PRZESUN}
    \frac{1}{4\pi(2n)^2} \int_{[0,2\pi)} \int_{(-n,n)^2}
     {\cal L}\left( [\Sigma \circ \tau_x \circ R_{\phi}]
     M_{\hat{\cal A}_{(-n,n)^2}^{[\alpha,\beta]}}\right) dx d \phi,
 \end{equation}
 where $\tau_x$ stands for the standard translation operator $\tau_x \mu(U) := \mu(U+x)$
 while $R_{\phi},\; \phi \in [0,2\pi)$ is the rotation by angle $\phi$ around $0$ and $\Sigma$
 is the reflection with respect to some fixed axis passing through the origin.
 By (\ref{OGG12}) it follows that
 $$ {\Bbb E} M_{n}^{[\alpha,\beta]}(U) < \infty $$
 for all bounded $U \subseteq {\Bbb R}^2.$ Applying Corollary A2.6.V in
 Daley \& Vere-Jones (1988) we conclude that the sequence
 of random measures $(M_{n}^{[\alpha,\beta]})_{n=1}^{\infty}$ is uniformly tight in
 $G_{{\Bbb R}^2}$ and, consequently, it contains a subsequence converging in
 law to some $M_{\infty}$ corresponding to a whole-plane polygonal field
 $\hat{\cal A}^{[\alpha,\beta]}_{\infty}.$ In view of (\ref{PRZESUN}) it
 is clear that $${\cal L}(\hat{\cal A}^{[\alpha,\beta]}_{\infty}) \in
 {\cal G}_{\tau}(\hat{\cal A}^{[\alpha,\beta]})$$ which completes
 the proof of the theorem. $\Box$



\subsection{Proof of Lemma \ref{NICHOLLS}}
 By the definition (\ref{THETAB}) of the $\beta$-tilted contour
 measure $\Theta^{[\beta]}$ it is enough to establish the assertion
 of the lemma for the henceforth assumed case $\beta := 2.$
 In order to establish (\ref{FKGREENA}) define the continuous-time
 random walk $Z_t$ in ${\Bbb R}^2$ with the following transition mechanism
 \begin{itemize}
   \item between critical events specified below move in a constant direction with speed $1,$
   \item with intensity given by $4$ times the covered length element update
         the movement direction, choosing the angle $\phi \in (0,2\pi)$ between
         the old and new direction according to the density $|\sin(\phi)| \slash 4,$
 \end{itemize}
 We start the random walk $Z_t$ at a given point $x$ and with a given
 initial velocity vector. Moreover, we choose the {\it loop-closing
 angle} $\phi^* \in (0,2\pi)$ according to the density $|\sin(\phi)|
 \slash 4$ and we draw an infinite {\it loop-closing} half-line $l^*$
 starting at $x$ and forming the angle $\phi^*$ with the initial velocity
 vector. Let $\hat{Z}_t$ be the random walk $Z_t$ killed whenever hitting
 its past trajectory or the loop-closing line $l^*.$ The directed nature
 of the random walk trajectories as constructed above requires considering
 for each contour $\theta$ two oriented instances $\theta^{\rightarrow}$
 (clockwise) and $\theta^{\leftarrow}$ (anti-clockwise). We claim that 
 for $x \in {\Bbb R}^2$ and $\theta \in {\cal C}$ with $x \in \Ver(\theta)$
 we have 
  $$
    8\pi dx \e^{-4 \lgth(e^*)} {\mathbb P}\left( \hat{Z}_t \mbox{ reaches $l^*$ and 
    the resulting contour falls into } d \theta^{\rightarrow} \right) $$
  \begin{equation}\label{RWR} 
   = \Theta^{[2]}(d\theta),
  \end{equation}  
  where $e^*$ stands for the last segment of $\theta^{\rightarrow}$
  counting from $x$ as the initial vertex, which is to coincide 
  with the segment of the loop-closing line $l^*$ joining its intersection
  point with $\hat{Z}_t$ to $x.$
  Clearly, the same relation holds then for $\theta^{\leftarrow},$
  hence adding versions of (\ref{RWR}) for $\theta^{\rightarrow}$
  and $\theta^{\leftarrow},$ which amounts to taking into account 
  two possible directions in which the random walk can move along 
  $\theta,$ will yield $2\Theta^{[\beta]}(d\theta)$ on the RHS.
  The relation (\ref{FKGREENA}) will easily follow by using the trivial upper
  bound $1$ for the probability on the LHS of (\ref{RWR}).

  To establish (\ref{RWR}), we observe that the probability element
  $${\mathbb P}\left( \hat{Z}_t \mbox{ reaches $l^*$ and 
  the resulting contour falls into } d \theta^{\rightarrow} \right)$$
  is exactly 
 \begin{equation}\label{RWR2}
   \frac{1}{4 [\mu \times \mu](\{(l,l^*)\;|\; l \cap l^* \in dx \})}
   \exp(-4\lgth(\theta \setminus e^*)) \prod_{i=1}^k d\mu(l[e_k]),
 \end{equation}
 where $e_1,\ldots,e_k$ are all segments of $\theta$ including $e^*,$  
 while $l[e_i]$ stands for the straight line determined by $e_i.$ Indeed,
 \begin{itemize}
  \item the prefactor $[4[\mu \times \mu](\{(l,l^*)\;|\; l \cap l^* \in dx \})]^{-1}$
        comes from the choice of the lines containing respectively the initial segment
        of $\theta^{\rightarrow}$ [counting from $x$] and $l^*,$ as well as from 
        the choice between two equiprobable directions on each of these lines,
  \item for the remaining segments we use the fact that, for any given straight
        line $l_0,$ $\mu(\{l\;|\; l \cap l_0 \in d\ell,\; \angle(l,l_0) \in d\phi\})
        = |\sin \phi| d\ell d\phi$ with $d\ell$ standing for the length element 
        on $l_0$ and with $\angle(l_0,l)$ denoting the angle between $l$ and
        $l_0,$ see Proposition 3.1 in Arak \& Surgailis (1989) as well as the argument
        justifying the dynamic representation in Section 4 ibidem. Note that the direction
        update intensity was set to $4$ to coincide with $\int_0^{2\pi} |\sin\phi| = 4.$ 
  \end{itemize}
  To get the required relation (\ref{RWR}) it is now enough to use
  (\ref{RWR2}), recall the definition of $\Theta$ and observe that 
  $[\mu \times \mu](\{(l,l^*)\;|\; l \cap l^* \in dx \}) = 2\pi dx$ as
  follows by standard integral geometry. This completes the proof 
  of (\ref{FKGREENA}).

 To proceed, let $\tilde{Z}_t$ be the random walk $Z_t$ killed whenever
 hitting its past trajectory, but not when hitting the loop-closing
 half-line $l^*.$ Define
 \begin{equation}\label{OKRSI}
  \varepsilon := - \lim_{T\to\infty} \frac{1}{T}
  \log {\Bbb P}\left( \tilde{\tau} > T \right),
 \end{equation}
 where $\tilde{\tau}$ is the lifetime of $\tilde{Z}_t$ or,
 in other words, the first moment when $Z_t$ hits its past
 trajectory. The existence of the limit in (\ref{OKRSI}) follows
 by a standard superadditivity argument, see Section 1.2 in Madras \&
 Slade (1993), and in fact $\varepsilon$ can be
 regarded as the connective constant for the self-avoiding version of
 the random walk $Z_t,$ see ibidem. It is easily checked that $\varepsilon > 0$
 since during each unit time of its evolution the walk $Z_t$ has a certain
 positive probability of hitting its past trajectory, uniformly bounded
 away from $0$ through time. To establish (\ref{ZSIGMA}) observe that,
 as in the argument above,
 $$ \Theta^{[2]}(\{ \theta\;|\; dx \in \Ver(\theta),\; \lgth(\theta) > R \}) \leq $$
 \begin{equation}\label{OGRAZSI}
     4 \pi dx {\Bbb P}\left(\tilde{Z}_t  \mbox{ survives up to time } R \slash 2 \right)
     \leq 4\pi dx {\Bbb P}(\tilde{\tau} > R \slash 2).
 \end{equation}
 The required relation (\ref{ZSIGMA}) follows now by (\ref{OGRAZSI}), (\ref{OKRSI})
 and by the observation that
 $$ \Theta^{[2]}(\{ \theta\;|\; {\bf 0} \in \innt \theta,\; \lgth(\theta) > R \}) \leq $$
 $$ \sum_{k=0}^{\infty} \Theta^{[2]}(\{ \theta\;|\; \Ver(\theta) \cap
    [B({\bf 0},k+1) \setminus B({\bf 0},k)] \neq \emptyset,\; \lgth(\theta) > \max(R,k) \}). $$
 The proof is complete. $\Box$

\paragraph{Acknowledgements}
 The author gratefully acknowledges the support of the Foundation for Polish
 Science (FNP).

\paragraph{References}
 \begin{description}
  \item{\sc Arak, T.} (1982) On Markovian random fields with finite number of values,
       {\it 4th USSR-Japan symposium on probability theory and mathematical statistics,
           Abstracts of Communications}, Tbilisi.
  \item{\sc Arak, T., Surgailis, D.} (1989) Markov Fields with Polygonal Realisations,
       {\it Probab. Th. Rel. Fields} {\bf 80}, 543-579.
  \item{\sc Arak, T., Surgailis, D.} (1991) Consistent polygonal fields,
       {\it Probab. Th. Rel. Fields} {\bf 89}, 319-346.
  \item{\sc Arak, T., Clifford, P., Surgailis, D.} (1993) Point-based polygonal models
        for random graphs, {\it Adv. Appl. Probab.} {\bf 25}, 348-372.
  \item{\sc Clifford, P., Nicholls, G.} (1994) A Metropolis sampler
        for polygonal image reconstruction, {\it available at}:\\
        {\tt http://www.stats.ox.ac.uk/~clifford/papers/met\_poly.html},
  \item{\sc Daley, D.J., Vere-Jones, D.} (1988) {\it An Introduction to the Theory of Point
        Processes}, Springer Series in Statistics, Springer-Verlag, New York.
  \item{\sc Fern\'andez, R., Ferrari, P., Garcia, N.} (1998)
        Measures on contour, polymer or animal models. A probabilistic
        approach. {\it Markov Processes and Related Fields} {\bf 4}, 479-497.
  \item{\sc Fern\'andez, R., Ferrari, P., Garcia, N.} (2002) Perfect simulation for
        interacting point processes, loss networks and Ising
        models. {\it Stoch. Proc. Appl.} {\bf 102}, 63-88.
  \item{\sc Nicholls, G.K.} (2001) Spontaneous magnetisation in the plane,
        {\it Journal of Statistical Physics}, {\bf 102}, 1229-1251.
  \item{\sc Schreiber, T.} (2004) Mixing properties for polygonal Markov fields in the plane,
        {\it submitted},
  \item{\sc G. R. Shorack, J. A. Wellner} (1986), {\it Empirical
        Processes with Applications to Statistics}, Wiley, New York.
  \item{\sc Surgailis, D.} (1991) Thermodynamic limit of polygonal models, {\it Acta applicandae
        mathematicae}, {\bf 22}, 77-102.
 \end{description}

\end{document}